\documentclass[12pt]{article}

\usepackage{amssymb}
\usepackage{amsmath}
\usepackage{graphicx}

\usepackage[usenames,dvipsnames]{color}

\def\rrd{{\rr^d}}

\def\1{^{-1}}
\def\one{{\bf1}}

\def\cald{{\mathcal{D}}}
\def\call{{\mathcal{L}}}

\def\calf{{\mathcal{F}}}
\def\calm{{\mathcal{M}}}
\def\calo{{\mathcal{O}}}
\def\calp{{\mathcal{P}}}

\def\calx{{\mathcal{X}}}

\def\rr{{\mathbb{R}}}
\def\nn{{\mathbb{N}}}
\def\9{{\infty}}
\def\lbb{{\lambda}}

\def\wt{\widetilde}
\def\ov{\overline}

\def\vf{{\varphi}}
\def\oo{{\omega}}

\def\pp{{\partial}}
\def\D{{\Delta}}
\def\vp{{\varepsilon}}

\def\barr{\begin{array}}
\def\earr{\end{array}}

\def\dd{\displaystyle}
\def\bk{\bigskip }

\def\n{\noindent }

\def\vsp{\vspace*{2mm}\\ }

\def\ff{\forall }

\def\({\left(}
\def\){\right)}
\def\<{\left<}
\def\>{\right>}

\newtheorem{theorem}{Theorem}[section]
\newtheorem{proposition}[theorem]{Proposition}
\newtheorem{lemma}[theorem]{Lemma}
\newtheorem{corollary}[theorem]{Corollary}

\newtheorem{remark}[theorem]{Remark}

\title{The evolution to equilibrium of solutions to~nonlinear Fokker-Planck equation}
\author{Viorel Barbu\thanks{Octav Mayer Institute of Mathematics of the Romanian Academy, Ia\c si, Romania} \and Michael R\"ockner\thanks{Fakult\"at f\"ur Mathematik, Universit\"at Bielefeld, D-33501 Bielefeld, Germany}\ \thanks{Academy of Mathematics and System Sciences, CAS, Beijing}}
\date{}

\begin{document}

\maketitle

\begin{abstract} \noindent  One proves   the $H$-theorem for mild solutions  to a nondegenerate, nonlinear Fokker-Planck equation
$$
u_t-\Delta\beta(u)+{\rm div}(E(x)b(u)u)=0, \ t\ge0, \ x\in\rr^d,\eqno(1)$$
and under appropriate hypotheses on $\beta,$ $E$ and $b$   the convergence in $L^1_{\rm loc}(\rr^d)$, $L^1(\rr^d)$, respectively,   for some $t_n\to\9$   of the solution $u(t_n)$ to an equi\-li\-brium state of the  equation for a large set of nonnegative  initial data in $L^1$. These results are new in the literature on nonlinear Fokker-Planck equations arising in the mean field theory and are also relevant to the theory of stochastic differential equations. As a matter of fact, by the above convergence result, it follows that  the solution to the McKean-Vlasov stochastic differential equation corresponding to (1), which is a {\it nonlinear distorted Brownian motion}, has  this equi\-li\-brium state as its unique invariant measure.\smallskip
\end{abstract}

\section{Introduction}

We shall study here the asymptotic behaviour of solutions $u=u(t,x)$ to the nonlinear Fokker-Planck equation

\begin{equation}\label{e1.1}
\barr{c}
u_t-\Delta\beta(u)+{\rm div}(Eb(u)u)=0\mbox{ in }(0,\9)\times\rr^d,\\
u(0,x)=u_0(x),\ x\in\rr^d,\earr  \end{equation}
 under the following hypotheses on the functions $\beta:\rr\to\rr,$ $E:\rr^d\to\rr^d$ and $b:\rr\to\rr$, where $1\le d<\9.$
\begin{itemize}
\item[(i)] $\beta\in C^1(\rr),\ \beta(0)=0,\ \gamma\le\beta'(r)\le\gamma_1,\
 \ff r\in\rr,$ for $0<\gamma<\gamma_1<\9.  $
	\item[(ii)] $b\in C_b(\rr)\cap C^1(\rr)$.
	\item[(iii)] $E\in L^\9(\rr^d;\rr^d)\cap W^{1,1}_{\rm loc}(\rrd;\rrd)$ and ${\rm div}\,E\in  (L^2(\rrd)+L^\9(\rrd))$.	
	\item[(iv)] {\it $E=-\nabla\Phi$, where $\Phi\in C(\rr^d)\cap W^{2,1}_{\rm loc}(\rr^d),\ \Phi\ge1,$ $\lim\limits_{|x|\to\9}\Phi(x)=+\9$ and there exists $m\in[2,\9)$ such that $\Phi^{-m}\in L^1(\rr^d)$.}
\end{itemize}Hypothesis (iv) means that system \eqref{e1.1} is conservative.

A typical example  is
$\Phi(x)=C(1+|x|^2)^\alpha,$ $ x\in\rr^d,$ with $\alpha\in\left(0,\frac12\right],$ for which we even have that ${\rm div}\,E\in L^\9$.

If (i)-(iv) hold, we prove the existence of solutions given by a nonlinear semi\-group $S(t)$, $t>0$, of contractions in $L^1(\rr^d)$ (Theorem \ref{t4.1}), which is po\-si\-ti\-vity  and mass preserving. If, (i)-(iv) and also (v) hold, where
\begin{itemize}
\item[(v)] $b(r)\ge b_0>0$ for $r\ge0,$
\end{itemize}
we prove the convergence of the solutions to equilibrium in $L^1_{\rm loc}(\rr^d)$, while (see Theorem \ref{t6.1}) the convergence in $L^1(\rr^d)$ is proved  if, in addition to (i)-(v), the following condition holds
 \begin{equation}\label{e1.2a}
{\rm(vi)}\hspace*{25mm}\gamma_1\Delta\Phi(x)-b_0|\nabla\Phi(x)|^2\le0, \ \mbox{ for a.e. } x\in\rr^d.\hspace*{10mm}
\end{equation}
An example of such a function $\Phi$ for $d\ge2$~is
\begin{equation}
\label{1.3}\Phi(x)=\left\{\barr{ll}
|x|^2\log|x|+\mu&\mbox{ for }|x|\le\delta,\vsp
\vf(|x|)+\eta|x|+\mu&\mbox{ for }|x|>\delta,\earr\right.
\end{equation} $\delta=\exp\(-\frac{d+2}{2d}\),$   and
\begin{equation}\label{e1.3prim}
\vf(r)=\delta^2\log\delta-\eta\delta+\dd\int^r_\delta h(s)ds, \end{equation}
for $r\ge\delta$, where $\mu,\eta>0$ are sufficiently large and $h$ is given by formula (A.8) in the Appendix to which we refer for more details.\newpage

 Equation \eqref{e1.1}, where $u$ is a probability density, is known in the literature as the nonlinear Fokker-Planck equation (NFPE) and it is relevant in the kinetic theory of statistical mechanics as a generalized mean field Smoluchowski equation for the case where the diffusion and transport coefficients depend on the density $u$. (See \cite{6}, \cite{7}-\cite{8} \cite{10}.) The case of the classical Smoluchowski equation is recovered for $b\equiv1$ and $\beta(r)\equiv r.$ In the case where the first order part in \eqref{e1.1} is given by a vector field independent of the spatial variable $x$, the existence and uniqueness of a kinetic, respectively generalized entropic, solution to \eqref{e1.1} in $L^1(\rr^d)$   was  proved  in \cite{8a}. In this paper, we   give an existence and uniqueness result for \eqref{e1.1}   in the sense of mild solutions in $L^1(\rr^d)$, i.e., given as a nonlinear semigroup $S(t)$, $t>0$,  in  $L^1(\rr^d)$ (see Proposition \ref{p2.2}). Its proof is dif\-fe\-rent from that in \cite{8a} and,  though it has an intrinsic interest in itself, it is   used subsequently to prove our main result about convergence to equilibrium and existence of a unique stationary solution to \eqref{e1.1}.  In \cite{4} (see, also, \cite{2}, \cite{3}),  a more general NFPE of the form
\begin{equation}\label{e1.2}
u_t-\sum^d_{i,j=1} D^2_{ij}(a_{ij}(x,u)u)+{\rm div}(b(x,u)u)=0\end{equation}
was studied under appropriate assumptions on $a_{ij}:\rr^d\times\rr\to\rr$ and \mbox{$b:\rr^d\times\rr\to\rr^d$.}
In the latter case, it is shown that, if $u_0$ is a probability density, the distributional mild solution $u$ to \eqref{e1.2} is  the probability density   of the law $\call_{X(t)}$ of the (probabilistically) weak solution to the McKean-Vlasov  stochastic differential equation (SDE)
\begin{equation}\label{e1.3}
dX(t)=b(X(t),u(t,X(t)))dt+\sqrt{2}\,\sigma(X(t),u(t,X(t)))dW(t),\end{equation}
where $\sigma\sigma^\bot=\frac12\,(a_{ij})^d_{i,j=1}$ and $X(0)$ has law $u_0dx$, where $dx=$ the Lebesgue measure on $\rr^d.$

In the special case   \eqref{e1.1}, SDE \eqref{e1.3} reduces to
\begin{equation}\label{e1.4a}
dX(t)=E(X(t))b(u(t,X(t)))dt+\frac1{\sqrt{2}}\,\(\frac{\beta(u(t,X(t)))}{u(t,X(t))}\)^{\frac12}dW(t),\end{equation}which, since $E=-\nabla\Phi$, is a nonlinear analogue of the SDE for the classical distorted Brownian motion, where $\beta=id$ and $b\equiv const.$ Hence, its solution $X(t)$, $t\ge0$, can be considered as a nonlinear distorted Brownian motion.

One of our motivations is to apply our asymptotic results {\it to find an inva\-riant $($probability$)$ measure for the nonlinear distorted Brownian motion on~$\rr^d$.}  So, Theorems \ref{t6.1} and \ref{t6.2} solve this problem and this is one of the main contributions of this work. Condition (vi) requires a certain balance between the strength of the (in general nonlinear) diffusion coefficient $\beta'$ and the strength of the nonlinear drift coefficient $b$ in terms of the {\it potential} $\Phi$. Without the additional condition (vi), there is in general no equilibrium on $L^1(\rr^d)$ for equation \eqref{e1.1}. Just consider the linear case $\beta=id$ and $E\equiv0$, so the case where \eqref{e1.1} is the heat equation. Hence, as in the linear case, we need a  big enough {\it negative} drift. Condition (vi) is, however, not optimal, because for the Fokker-Planck equation associated to the classical Ornstein-Uhlenbeck process on $\rr^d$, it does not hold, though the standard Gaussian measure is its  equilibrium measure.

We would like to mention here another special case of \eqref{e1.1}, namely with $\beta(u)=u^m,$ $m>1,$ $b\equiv const.$ and $E(x)=x$, which is not covered by our results, but was deeply analyzed in \cite{9prim}. In this case, the equilibrium is given through an explicit formula and the decay rate in $L^1$-distance is calculated in \cite{9prim}. So, the approach is completely different from ours  which is to prove the so-called $H$-theorem (see below) to show convergence of solutions to a unique equilibrium of \eqref{e1.1} in $L^1(\rr^d)$ as $t\to\9.$   A general result combining \cite{9prim}, the linear case and ours including convergence rates is still to be proved and will be subject to our future study.
 As explained in detail in \cite[Section 2]{4}, the nonlinear Fokker-Planck equation \eqref{e1.1} is a (very singular) special case (called {\it Nemytskii type}) of a general nonlinear Fokker-Planck-Kolmogorov equation in the sense of Section 6.7(iii) in \cite{0} and of  \cite{18prim}, \cite{18second}, where the solutions are measure-valued and the coefficients depend on these solutions. There is a
number of papers where existence of and convergence to equilibria are studied (see, e.g., \cite{0prim} and \cite{14primm} and the references therein). However, in these papers the dependence of the coefficients on the measures is assumed to be linear or Lipschitz continuous in weighted variation norm, which is never fulfilled in our Nemytskii-type case. So, these results do not apply here.

 The main objective of this work is to study the asymptotic behaviour of a~solution $t\to u(t)$ for $t\to\9$ and prove the so called  \mbox{$H$-theorem} for~the  NFPE \eqref{e1.1}, that is, prove the existence of a Lyapunov function  \mbox{$V:D(V)\subset L^1_{\rm loc}(\rr^d)\to\rr$} for \eqref{e1.1} and prove, for a certain  class of $u_0\in L^1,$ $u_0\ge0$,  the $\oo$-limit set\vspace*{-2mm}
\begin{equation}\label{e1.4}
\oo(u_0)=\left\{w=\lim_{n\to\9}u(t_n)\mbox{ in }L^1_{\rm loc} (\rr^d),\ \{t_n\}\to\9\right\}\end{equation}
is nonempty.  This is proved in Sections \ref{s4} and \ref{s5}  under assumptions \mbox{(i)-(v).}

Moreover, if (vi) also holds, we shall prove in Section \ref{s6}  that,  for $u_0\in\calm\cap\calp$ (see \eqref{e2.2}, \eqref{e2.24a}), the orbit  $\{u(t);\ t\ge0\}$ is compact in $L^1$ and so the corresponding $\oo$-limit set $\wt\oo(u_0)=\left\{w=\lim\limits_{n\to\9}u(t_n)\mbox{ in }L^1,\ \{u_n\}\to\9\right\}$ is nonempty and   reduces to a single element $u_\9$, which is a stationary solution to \eqref{e1.1}.
Furthermore, $u_\9$ is a probability density, if so is $u_0$. As a consequence,   $u_\9dx$ is an invariant measure for SDE \eqref{e1.4a}, i.e., if $u_0=u_\9$, then {\it the nonlinear distorted Brownian motion $X(t)$, $t\ge0$, has the law $u_\9dx$, $\ff t\ge0.$}

The $H$-theorem amounts to saying that the function
\begin{equation}\label{e1.5}
V(u)=-\wt S[u]+F[u],\ u\in L^1(\rr^d),\end{equation}where $\wt S$ is the entropy of the system and $F$ is the   mean field energy, is a Lyapunov function for \eqref{e1.1}, that is, monotonically decreasing in time on the solutions to \eqref{e1.1}. In our case,
\begin{equation}\label{e1.8}
\wt S[u]=\int_{\rr^d}\eta(u(x))dx,\ F(u)=\int_{\rr^d}\Phi(x)u(x)dx,\end{equation}
where
 $\eta(r)=-\int^r_0d\tau\int^1_\tau
\frac{\beta'(s)}{sb(s)}\ ds,\ r\ge0.$

This form of the Lyapunov theorem comes from the classical $H$-theorem and is consistent with the Boltzmann thermodynamics (see, e.g., \cite{6}, \cite{7}, \cite{10}), in which case $\beta'\equiv b\equiv const.$, so $\wt S$ in \eqref{e1.8} reduces to the classical Boltzmann-Gibbs entropy.
In the literature  on NFPE arising in the mean field theory,   the $H$-theorem is often invoked, but in most cases  its proof is   formal because, in general, the NFPE \eqref{e1.1} has not a classical solution and so the computation is not ri\-go\-rous. By our knowledge,  this paper contains  the first ri\-go\-rous mathematical result on the $H$-theorem for NFPE.

In fact, here the basic functional space for the well-posedness is $L^1(\rr^d)$ and, in general, the space of the maximal spatial regularity for $u$ is the Sobolev space $W^{1,q}(\rr^d)$, $1<q\le\frac d{d-2}$, (which happens in the special case of the porous media equation $b\equiv0$, $a_{ij}(u)u\equiv \delta_{ij}\beta(u)).$ This low regularity precludes the classical argument involving regular Lyapunov functions. However, the situation is different for linear FPE where, in the last decades, many convergence results to equili\-brium  were obtained. We refer to the monographs \cite{1a}, \cite{22} and, e.g., to \cite{15a}, \cite{9prim}, \cite{16a}, \cite{25a}, as well as the references therein.

Here, the convergence of $S(t)u_0$ for $t\to\9$ to an equilibrium state is proved under nondegenerancy assumption (i) for $\beta$. In the degerate case, $\beta'>0$ on $[0,\9)$, one expect, however, that the omega limit set $\wt\oo(u_0)$ is nonempty and is a compact attractor for $S(t)$.   (We refer to   \cite{26a} for a theory of infinite dimensional attractor.)

Let us now explain the structure of the paper. The first part is concerned with the well-posedness of NFPE \eqref{e1.1} in $L^1(\rr^d)$ via the theory of non\-linear semigroups of contractions in $L^1(\rr^d)$,   i.e., the construction of such a semigroup $S(t)$, $t>0$, so that $t\mapsto S(t)u_0$   a continuous function \mbox{$u:[0,\9)\to L^1(\rr^d)$} given as the limit of the finite difference scheme associated with \eqref{e1.1} (the so called {\it mild} solution). Moreover,   $u$ is obtained as  the limit in $L^1(\rr^d)$ of the smooth solutions $\{u_\vp\}_{\vp>0}$ to an approximating equation associated with \eqref{e1.1}. The correspon\-ding  result given in Proposition \ref{p2.1}  is not essentially new since, as mentioned earlier, a similar existence result   was previously esta\-blished in \cite{2}-\cite{6a}, \cite{8a}.  However, we have developed here a semigroup approach to NFPE \eqref{e1.1} necessary for the treatment of   the asymptotic behaviour of solutions.   In fact, in the second part of the work we shall prove  under assumptions (i)-(v)   the $H$-theorem for \eqref{e1.1}  (Theo\-rem \ref{t4.1}). The $\oo$-limit set is a singleton $\{u_\9\}$  and the   invariant measure of the solution $X(t)$, $t\ge0$, of SDE \eqref{e1.4a}  if, additionally, the balance condition (vi) holds  (Theorem \ref{t6.1}). A main point to prove the latter is to show that $S(t)$  is also a contraction on the weighted $L^1$ space with the potential $\Phi$ from condition (iv) as its weight (see Lemma \ref{l6.2}).

Finally, we prove that the equilibrium $u_\9$ from
Theorem \ref{t6.1} is indeed the unique solution of the stationary version of \eqref{e1.1} in the sense of distributions (Theorem \ref{t6.2}) and, as a consequence, that the stationary nonlinear distorted Brownian motion is unique in law (Theorem \ref{t6.4}).

\medskip\noindent{\bf Notation.}  For $p\in[1,\9)$, $L^p(\rr^d)$ - simply denoted $L^p$, is the space of all Lebesgue $p$-summable functions on $\rr^d$. The norm in $L^p$ is denoted by~$|\cdot|_p$. Similarly, if $\calo$ is a Lebesgue measurable set, $L^p(\calo)$ is the space of all $p$-summable functions on $\calo$. By $L^p_{\rm loc}(\rr^d)$ we denote the space of Lebesgue measurable functions $u:\rr^d\to\rr$ which are in $L^p(\calo)$ for every bounded measurable subset $\calo\subset\rr^d$. ($L^p_{\rm loc}$ is endowed with a standard locally convex metrizable topology.) The scalar product of $L^2$ is denoted by $\<\cdot,\cdot\>_2$. If $\calo$ is an open subset of $\rr^d$, we denote by $\cald'(\calo)$ the space of Schwartz distributions on $\calo$ and by $W^{1,p}(\calo)$ the Sobolev space $\{u\in L^p(\calo),\ D_iu\in L^p(\calo)$ for $i=1,...d\}$, where $D_i=\frac\pp{\pp x_i}$ is taken in the sense of Schwartz distributions. We set also $H^k(\calo)=W^{k,2}(\calo)$, $k\in\nn$. We denote the Euclidean norm of $\rr^d$ by $|\cdot|$, if there is no possible confusion, and by $C_b(\rr)$ and $C_b(\rr^d,\rr^d)$ the spaces of continuous and bounded functions from $\rr$ to itself and, respectively, from $\rr^d$ to $\rr^d$. By $C^1(\rr)$ we denote the space of continuously differentiable real valued functions.

\section{Existence of mild solutions for NFPE \eqref{e1.1}}\label{s2}
\setcounter{equation}{0}

 Consider in the space $L^1=L^1(\rr^d)$ the operator $A_0:D(A_0)\subset L^1\to L^1$, defined~by
\begin{equation}\label{e2.1}
\barr{rcl}
A_0u&=&-\D\beta(u)+{\rm div}(Eb(u)u),\ \ff u\in D(A_0),\vsp
D(A_0)&=&\{u\in L^1;\ -\D \beta(u)+{\rm div}(Eb(u)u)\in L^1\}.\earr\end{equation}Here, the differential operators $\D$ and $\rm div$ are taken in the sense of Schwartz distributions, i.e., in $\cald'(\rr^d)$. Obviously, the operator $(A_0,D(A_0))$ is closed on~$L^1$.

By Hypotheses (i)-(iii), we see that $\beta(u),Eub(u)\in L^1$, $\ff u\in L^1$, and so $-\D\beta(u),{\rm div}(Eub(u))\in\cald'(\rr^d)$ for all $u\in  L^1$.

\begin{proposition}\label{p2.1} Assume that Hypotheses {\rm(i)-(iv)} hold. Then,
\begin{eqnarray}
&R(I+\lbb A_0)=L^1,\ \ff\lbb>0,\label{e2.2}\end{eqnarray}
and there is an operator $J_\lbb:L^1\to L^1$ such that $J_\lbb(0)=0,\ \lbb>0$, and
\begin{eqnarray}
&J_{\lbb_2}(f)=J_{\lbb_1}\(\dd\frac{\lbb_1}{\lbb_2}\ f+\(1-\dd\frac{\lbb_1}{\lbb_2}\)J_{\lbb_2}(f)\),\ \ff\lbb_1,\lbb_2>0,\label{e2.3a}\\[1mm]
&(I+\lbb A_0)J_\lbb(f)=f,\ \ \ff f\in  L^1,\ \lbb>0,\label{e2.3aa}\\[1mm]
&|J_\lbb(f_1)-J_\lbb(f_2)|_1\le|f_1-f_2|_1,\ \ff\lbb>0,\ f_1,f_2\in L^1.\label{e2.3}\end{eqnarray}
 Furthermore,
 \begin{equation}
\overline{D(A)}=L^1,\label{e2.3d}
 \end{equation}where $\overline{\raise 5pt\hbox{\ \ \ \ }}$ denotes the closure in $L^1$ and $A$ is the operator defined by formula \eqref{e2.6a} below.
Moreover,
\begin{eqnarray}\label{e2.4}
&\dd\int_{\rr^d}J_\lbb (f)dx=\int_{\rr^d} f(x)dx,\ \ff  f\in L^1,\\[2mm]
 \label{e2.5}
&\dd J_\lbb (f)\ge0,\ \mbox{ a.e. in }\rr^d\mbox{ if } f\ge0,\mbox{ a.e.  in }\rr^d.\label{e2.6c}\end{eqnarray}
\end{proposition}

The proof of Proposition \ref{p2.1} will be given in Section \ref{s3}.

We note that $J_{\lbb_1}(L^1)=J_{\lbb_2}(L^1)$, $ \ff\lbb_1,\lbb_2>0.$ We are led to introduce the operator   $A:D(A)\subset L^1\to L^1$,
\begin{equation}
\label{e2.6a}
Au=A_0u,\ \ff u\in D(A)=J_{\lbb_0}(L^1),\ \ff\lbb>0,
\end{equation}where $\lbb_0>0$ is arbitrary. {Hence, $D(A)\subset D(A_0)$ and} taking into account \eqref{e2.3a}, it follows that $D(A)$ is independent of $\lbb_0$.

By \eqref{e2.2}-\eqref{e2.3d}, it follows that {\it $A$ is $m$-accretive in $L^1$.}   This means (see, e.g. \cite{15a}, p.~97) that
$|u-v+\lbb(Au-Av)|_1\ge|u-v|_1,\ \ff u,v\in D(A),$ $ \lbb>0,$ and $R(I+\lbb A)=L^1$, $\ff\lbb>0$ (equivalently, for some $\lbb>0$). We have
\begin{equation}
\label{e2.6aa}
(I+\lbb A)\1u=J_\lbb(u),\ \ff u\in L^1,\ \lbb>0.
\end{equation}We note that $A$ is an accretive section of $A_0$ and if $(I+\lbb A_0)\1$ is single valued, then $A=A_0$. As shown in \cite{8aa} (Proposition 2.4), this happens for instance if, besides (i)--(iii), the following conditions hold
\begin{equation}
\label{e2.10a}
{\rm div}\,E\in L^{m}_{\rm loc},\ m>\frac d2,\ |rb'(r)+b(r)|\le\alpha\beta'(r),\ \ \ff\,r\in\rr;\ \alpha>0.
\end{equation}

Consider now the Cauchy problem associated with $A$, that is,
\begin{equation}\label{e2.6}
\barr{l}
\dd\frac{du}{dt}+Au=0,\  t\ge0,\vsp
u(0)=u_0.\earr\end{equation}A continuous function $u:[0,\9)\to L^1$ is said to be a {\it mild solution to equation \eqref{e2.6}} if
\begin{equation}\label{e2.7}
u(t)=\lim_{h\to0}u_h(t)\mbox{ in }L^1,\end{equation}uniformly on compacts of $[0,\9)$, where $u^{1}_h=u_0$, and
\begin{eqnarray}\label{e2.8}
& u_h(t)=u^i_h,\ t\in[ih,(i+1)h),\ i=0,1,..., \\[2mm]
 \label{e2.9}
& u^{i}_h+hAu^{i}_h=u^{i-1}_h,\ i=0,...\end{eqnarray}Since $A$ is $m$-accretive, we have by the Crandall \& Liggett theorem (see, e.g., \cite{1}, p. 141)  the following existence result for problem \eqref{e2.6}.

\begin{proposition}\label{p2.2} Under Hypotheses {\rm(i)-(iv)}, for every $u_0\in L^1(\rrd)$ there is a unique mild solution $u=S(t)u_0$ to  \eqref{e2.6}. Moreover,  one has\newpage	
	
\begin{equation}\label{e2.10}
u(t)=\lim_{n\to\9}\(I+\frac tn\,A\)^{-n}u_0,\ \  \ff\,t\ge0, \end{equation}uniformly on bounded intervals of $[0,\9)$ in the strong topology in $L^1$.
One   also has that
\begin{eqnarray}\label{e2.11}
&\dd\int_{\rr^d}u(t,x)dx=\int_{\rr^d}u_0(x)dx,\ \ff t\ge0,\\[2mm]
 \label{e2.12a}
&\dd u(t,x)\ge0,\mbox{ a.e. on }(0,\9)\times\rr^d\mbox{ if }u_0\ge0,\mbox{ a.e. in }\rr^d.\end{eqnarray}
\end{proposition}Taking into account that by \eqref{e2.6a}--\eqref{e2.6aa}, equation \eqref{e2.8} can be written as
\begin{equation}
\label{e2.14prim}u^i_h-h\Delta\beta(u^i_h)+h\ {\rm div}(Eb(u^i_h)u^i_h)=u^{i-1}_h\mbox{ in }\cald'(\rrd),
\end{equation}
the function $u$ will be called {\it mild solution} to NFPE \eqref{e1.1}.

In particular, it follows by \eqref{e2.11}, \eqref{e2.12a} that, for each $t\ge0$, $u(t,\cdot)$ is a probability density if so is $u_0$.

We note that \eqref{e2.11}-\eqref{e2.12a} follow by \eqref{e2.4}-\eqref{e2.5}  and \eqref{e2.10}.

The map $t\to S(t)u_0$ is a continuous semigroup of contractions on $L^1$, that is,
\begin{eqnarray}\label{e2.12}
&\dd S(t)u_0=u(t)=\lim_{n\to\9}\(I+\frac tn\,A\)^{-n}u_0,\ \ff t\ge0,\\[1mm]
 \label{e2.13}
&\dd S(t+s)u_0=S(t)S(s)u_0,\ \ff t,s\ge0,\ u_0\in L^1,\\[1mm]
 \label{e2.14}
&\dd\lim_{t\to0}S(t)u_0=u_0\mbox{ in }L^1,\\
 \label{e2.15}
&\dd|S(t)u_0-S(t)\bar u_0|_1\le|u_0-\bar u_0|_1,\ \ff t\ge0,\ u_0,\bar u_0\in L^1.\end{eqnarray}
If
\begin{equation}\label{e2.16}
\calp=\left\{u\in L^1;\ u\ge0,\mbox{ a.e. in }\rrd,\  \int_{\rr^d}u(x)dx=1\right\},\end{equation}
we see by \eqref{e2.11}-\eqref{e2.12} that
\begin{equation}\label{e2.17}
S(t)(\calp)\subset\calp,\ \ff t\ge0,\end{equation}and, since $J_\lbb(0)=0,$ that
\begin{equation}\label{e2.24'}
S(t)(0)=0,\ t\ge0.\end{equation}
Since, for every $i$ and $h$ the function $u^{i}_h\in D(A)$ is a solution to \eqref{e2.9} in the sense of  distributions, i.e. in the space $\cald'(\rr^d)$, it follows also that the mild solution $u$ to \eqref{e2.6} is a solution to NFPE \eqref{e1.1} in the sense of Schwartz distributions on $(0,\9)\times\rr^d$, that is,
\begin{equation}\label{e2.18}
\barr{l}
\dd\int^\9_0\!\!\int_{\rr^d}(u\vf_t+\beta(u)\D\vf+Eb(u)u\cdot\nabla\vf)dx\,dt\\
\qquad\qquad+\dd\int_\rrd u_0\vf(t,x)dx=0,\ \ff\vf\in\cald([0,\9)\times\rr^d),\earr\end{equation}where $\cald((0,\9)\times\rr^d)$ is the space of infinitely differentiable functions on \mbox{$(0,\9)\times\rr^d$} with compact support. 

It should be emphasized, however, that the solution $u$ to NFPE \eqref{e1.1} \mbox{exists} and is unique in the class of mild solutions corresponding to the ope\-ra\-tor $A$ and not in the space of Schwartz distributions on $(0,\9)\times\rr^d.$ In other words, it is dependent on $\{J_\lbb\}$ which in our case is the limit of $(I+\lbb (A_0)_\vp)\1$ in $L^1$, where $(A_0)_\vp$ is a smooth approximation of $A_0$. However, as $u=S(t)u_0$ is $L^1$-valued continuous, then, as shown in \cite{7a}, \cite{8ab} under the additional condition that $u_0\in L^\9$,
  it is unique in this case in the class of distributional solutions $u\in L^\9((0,\9)\times\rr^d)\cap L^1((0,\9)\times\rr^d)$   and so it is unique in the class of all mild solutions with $u_0\in L^1\cap L^\9$. The semigroup $S(t)$ can be viewed, therefore, as the Fokker--Planck flow generated by equation \eqref{e1.1} which is uniquely defined on the space $L^1\cap L^\9$.

We consider the following subspace of $L^1$
\begin{equation}\label{e2.24a}
\barr{rcl}
\calm&=&\left\{u\in L^1;\dd\int_{\rr^d}\Phi(x)|u(x)|dx<\9\right\}
\earr\end{equation}
with the norm
\begin{equation}\label{e2.21}
\|u\|=\dd\int_{\rr^d}\Phi(x)|u(x)|dx,\,\ff u\in\calm.
\end{equation}We also set
 $\ \calm_+=\{u_0\in\calm;\ u_0\ge0,\mbox{ a.e. on }\rrd\}.$

It turns out that the semigroup $S(t)$ leaves invariant $\calm$. More precisely, we prove in Section \ref{s3}:
\begin{proposition}\label{p2.3}  Assume that Hypotheses {\rm(i)-(iv)}  hold and that ${\rm div}\,E\in L^\9$.   Then
\begin{equation}\label{e2.21a}
\|S(t)u_0\|\le\|u_0\|+ \rho t|u_0|_1,\ \ff u_0\in\calm, \end{equation}where
$\rho=\gamma_1(m+1)|\Delta\Phi|_\9+|b|_\9(1+m)^2|E|^2_\9.$
\end{proposition}

\begin{remark}\label{r2.4}\rm Proposition  \ref{p2.3} remains valid if, in addition to Hypotheses (i)-(iii), we assume, instead of (iv),
\begin{itemize}
\item[(iv)$'$] $E_0=\sup\limits_{x\in\rr^d}\,|E(x)\cdot x|<\9,$
    \end{itemize}
    but we have to replace $\calm$ by
    $$\calm_2=\left\{u\in L^1:\|u\|_2=\int_{\rr^d}|x|^2|u(x)|dx<\9\right\}$$
   and we have to replace $\rho$ in Proposition \ref{p2.3} by $\wt\rho:=2(d\gamma_1+E_0|b|_\9) $ (see Remark \ref{r3.3} below). The assumption (iv), in particular that $E$ is the negative of the gradient of a positive function, becomes, however, important for Sections \ref{s4}-\ref{s6} below, i.e., to prove the $H$-Theorem.
\end{remark}

\section{Proof of Propositions \ref{p2.1} and \ref{p2.3}}\label{s3}
\setcounter{equation}{0}

As mentioned earlier, one can derive Proposition \ref{p2.1} from similar results established in \cite{3}, \cite{4}.
However, for later use we shall prove it by a constructive regularization technique already developed in the above works.    Namely, we define, for each $\vp>0$, the operator
$(A_0)_\vp:D((A_0)_\vp)\subset L^1\to L^1$,
\begin{eqnarray}
(A_0)_\vp u\!\!\!&=&\!\!\!-\D(\beta(u))+\vp\beta(u)+{\rm div}(E_\vp b^*_\vp(u)),\label{e3.1}\\[2mm]
D((A_0)_\vp)\!\!\!&=&\!\!\!\{u\in L^1,\ -\D(\beta(u))+\vp\beta(u)\label{e3.2}+{\rm div}(E_\vp   b^*_\vp(u))\in L^1\}.\qquad  \end{eqnarray}
Here $\D$ and div are taken in the sense of Schwartz distributions and
\begin{equation}\label{e3.2e}
 b_\vp\equiv   b*\rho_\vp,\ \
b^*_\vp(r)\equiv\dd\frac{b_\vp(r)r}{1+\vp|r|},\ r\in\rr,
\end{equation}
where  $\rho_\vp(r)\equiv\frac1\vp\ \rho\(\frac r\vp\)$, $\rho\in C^\9_0(\rr)$, $\rho\ge0$, is a  standard mollifier. Moreover,
$$E_\vp =-\nabla\Phi_\vp ,\ \ \Phi_\vp (x)\equiv\frac{\Phi(x)}{(1+\vp\Phi(x))^m}.$$Then $\Phi_\vp\in L^2$, since $m\ge2$, and
\begin{equation}\label{e33a}
E_\vp =E(1+\vp\Phi)^{-m}-m\vp\Phi E(1+\vp\Phi)^{-(m+1)}
\end{equation}and, therefore, by Hypothesis (iv),

\begin{equation}\label{e3.4x}
\barr{c}
E_\vp \in (L^\9\cap L^1)(\rr^d;\rr^d)\vsp
|E_\vp(x)|\le(1+m)|E(x)|,\   \
\dd\lim_{\vp\to0}E_\vp(x)=E(x),\ \mbox{ for a.e. } x\in\rr^d,\vsp
\vp^{m}|E_\vp|\le(1+m)|E|_\9\Phi^{-m},\ \ff\vp>0.\earr\end{equation}
 We also note that $b^*_\vp,b_\vp$ are bounded and Lipschitz   and that, for $\vp\to0$, \begin{equation}
b^*_\vp(r)\to b(r)r \ \mbox{ uniformly on compacts}. \label{e3.2a}\end{equation}
Obviously, the operator $((A_0)_\vp,D((A_0)_\vp))$ is closed on $L^1$.

\begin{lemma}\label{l3.1}Assume that Hypotheses {\rm(i)-(iv)} hold. Then
\begin{equation}\label{e3.6a}
R(I+\lbb(A_0)_\vp)=L^1,\ \ff\lbb>0,
\end{equation}and there is an operator $J^\vp_\lbb:L^1\to L^1$ such that $J^\vp_\lbb(0)=0$ and \eqref{e2.3a}--\eqref{e2.3} hold. Namely,
\begin{eqnarray}
&J^\vp_{\lbb_2}(f)=J^\vp_{\lbb_1}\(\dd\frac{\lbb_1}{\lbb_2}\,f+\(1-\frac{\lbb_1}{\lbb_2}\)J^\vp_{\lbb}(f)\),\ \ff\lbb_1,\lbb_2>0,\label{e3.6aa}\\[1mm]
&(I+\lbb(A_0)_\vp)J^\vp_{\lbb}(f)=f,\ \ff f\in L^1,\ \ff\vp>0,\label{e3.6aaa}\\[2mm]
	&|J^\vp_\lbb(f_1)-J^\vp_\lbb(f_2)|_1\le|f_1-f_2|_1,\ \ff f_1,f_2\in L^1,\ \lbb>0,\label{e3.6aaaa}\\[2mm]
	&J^\vp_\lbb(f)\ge0,\mbox{ a.e. in }\rrd\mbox{ if }f\ge0,\mbox{ a.e. in }\rrd,\ \ff\lbb\in(0,\lbb_1),\label{e3.6aaaaa}\\[1mm]
	&\dd\int_\rrd J^\vp_\lbb(f)dx=\int_\rrd f\,dx,\ \ff\lbb>0,\ \ff\,f\in L^1.\label{e3.6aaaaaa}
	\end{eqnarray} 
Moreover, there is $\lbb_0>0$ independent of $f\in L^1$ such that, for all $\lbb\in (0,\lbb_0)$,
\begin{equation}\label{e3.3}
\lim_{\vp\to0}J^\vp_\lbb(f)=J_\lbb(f)\mbox{\ \ in }L^1,\ \ff f\in L^1,\end{equation}where $J_\lbb$ satisfies \eqref{e2.3a}--\eqref{e2.3} and \eqref{e2.4}, \eqref{e2.5}.
\end{lemma}

As in the case of the operator $A$, we define (see \eqref{e2.6a})
\begin{equation}
\label{e3.7a}
A_\vp u=(A_0)_\vp u,\ \ff u\in D(A_\vp)=J^\vp_\lbb(L^1).
\end{equation}
Then, {\it Lemma {\rm\ref{l3.1}} implies that $A_\vp$ is $m$-accretive in $L^1$ and $(I+\lbb A_\vp)\1=J^\vp_\lbb$.} Moreover, by \eqref{e3.3} it follows that
\begin{equation}
\label{e3.7aa}
\lim_{\vp\to0}(I+\lbb A_\vp)\1 f=J_\lbb(f)\mbox{\ \ in }L^1,\ \ff f\in L^1,\ \mbox{for $\lbb\in(0,\lbb_0).$}
\end{equation}

\noindent{\bf Proof of Lemma \ref{l3.1}.} We  fix $f\in L^2\cap L^1$ and consider the equation $u+\lbb (A_0)_\vp u=f,$ that~is,\newpage
\begin{equation}\label{e3.7}
u-\lbb\D(\beta(u))+\vp\lbb\beta(u)+\lbb\ {\rm div}(E_\vp  b^*_\vp(u))=f\mbox{ in }\cald'(\rr^d).\end{equation}
To solve equation \eqref{e3.7}, we consider the equation
\begin{equation}\label{e3.8}
(\vp I-\D)^{-1}u+\lbb\beta(u)+\lbb(\vp I-\D)^{-1}{\rm div}(E_\vp  b^*_\vp(u)) =(\vp I-\D)^{-1}f\mbox{ in }L^2.\end{equation}Clearly, a solution of \eqref{e3.8} satisfies \eqref{e3.7} in $L^2$. We set
\begin{equation}\label{e3.9}
\barr{rcl}
F_\vp(u)&=& (\vp I-\D)^{-1}u,\ G(u)=\lbb\beta(u),\ u\in L^2,\\
G_\vp(u)&=&\lbb(\vp I-\D)^{-1}({\rm div}(E_\vp  b^*_\vp(u))),\ u\in L^2,\earr\end{equation}and note that $F_\vp$ and $G$ are accretive and continuous in $L^2.$

We also have by Hypotheses (ii)-(iii) that $G_\vp$ is continuous in $L^2$ and
\begin{equation}\label{e3.10}
\barr{l}
\dd\int_{\rr^d}(G_\vp(u)-G_\vp(\bar u))(u-\bar u)dx\vspace*{-2mm}\\
\qquad=-\lbb\dd\int_{\rr^d}E_\vp (b^*_\vp(u)-b^*_\vp(\bar u))\cdot\nabla(\vp I-\D)^{-1}(u-\bar u))dx\vsp
\qquad\dd\ge -C_\vp \lbb |u-\bar u|_2|\nabla(\vp I-\D)^{-1}(u-\bar u)|_2,\ \ff u,\bar u\in L^2(\rr^d),\earr\end{equation}for some positive constant $C_\vp=0\(\frac1\vp\)$.
 Moreover, we have
\begin{equation}\label{e3.11}
\int_{\rr^d}(\vp I-\D)^{-1}uu\,dx=\vp|(\vp I-\D)^{-1}u|^2_2+|\nabla(\vp I-\D)^{-1}u|^2_2,\ \ff u\in L^2.\end{equation}
By \eqref{e3.8}-\eqref{e3.11}, we see   that, for $u^*=u-\bar u$, we have
$$\barr{l}
(F_\vp(u^*)+G_\vp(u)-G_\vp(\bar u)+G(u)-G(\bar u),u^*)_2\vsp
\qquad
\ge\lbb\gamma|u^*|^2_2+|\nabla(\vp I-\D)^{-1}u^*|^2_2+\vp|(\vp I-\D)^{-1}u^*|^2_2\vsp
\qquad-C_\vp \lbb |u^*|_2|\nabla(\vp I-\D)^{-1}u^*|_2.\earr$$This implies that $F_\vp+G_\vp+G$ is  accretive and   coercive on $L^2$ for $\lbb<\lbb_\vp$, where $\lbb_\vp$ is sufficiently small. Since this operator is continuous and accretive, it follows that it is $m$-accretive and, therefore, surjective (because it is coercive). Hence, for each $f\in L^2\cap L^1$ and $\lbb<\lbb_\vp$, equation \eqref{e3.8} has a unique solution $u_\vp\in L^2$.
Since $u_\vp\in L^2$, $b^*_\vp(r)\le C_\vp|r|, $ $r\in\rr$, and $E_\vp \in L^\9$, by \eqref{e3.7} we see that $\beta(u_\vp)\in H^1(\rr^d)$,
  whence  by (i) we have
\begin{equation}\label{e3.12c}u_\vp\in H^1(\rr^d).\end{equation}\newpage
Multiplying  \eqref{e3.7} by  $u_\vp$ and $\beta(u_\vp)$, respectively, integrating over $\rr^d$ and using hypothesis (i) (part $\beta'\ge\gamma$), we get after some calculation that, for $\lbb<\lbb_1$ with $\lbb_1$ small enough, 
\begin{equation}\label{e3.13a}
|u_\vp|^2_2+\lbb|\nabla\beta(u_\vp)|^2_2+\lbb|\nabla u_\vp|^2_2+
\vp\lbb|\beta(u_\vp)|^2_2\le C_{\lbb_1}|f|^2_2,
\end{equation}where $C_{\lbb_1}$ is independent of $\vp$.

We denote by $u_\vp(f)\in H^1(\rr^d)$ the solution to \eqref{e3.8} for $f\in L^2\cap L^1$ and prove that
\begin{equation}\label{e3.13b}
|u_\vp(f_1)-u_\vp(f_2)|_1\le|f_1-f_2|_1,\ \ff f_1,f_2\in L^1\cap L^2.
\end{equation}Here is the argument. We set $u=u_\vp(f_1)-u_\vp(f_2),$ $f=f_1-f_2$.   By \eqref{e3.7}, we have, for $u_i=u_\vp(f_i)$, $i=1,2,$
\begin{equation}\label{e3.13bb}
\barr{r}
u-\lbb\Delta(\beta(u_1)-\beta(u_2))
+\vp\lbb(\beta(u_1)-\beta(u_2))\vsp
+\lbb\,{\rm div}(E_\vp (b^*_\vp(u_1)-b^*_\vp(u_2)))=f\ \mbox{ in }L^2.\earr\end{equation}
Proceeding as in \cite{4} (see, also, \cite{7c}), we consider   the Lipschitzian function $\calx_\delta:\rr\to\rr,$
\begin{equation}\label{e3.18az}
\calx_\delta(r)=\left\{\barr{rl}
1&\mbox{ for }r\ge\delta,\\\dd\frac r\delta&\mbox{ for }|r|<\delta,\\-1&\mbox{ for }r<-\delta,\earr\right.\end{equation}
where $\delta>0$.
 We set
  $$F_\vp=\lbb \nabla(\beta(u_1)-\beta(u_2))
  -\lbb E_\vp  (b^*_\vp(u_1)-b^*_\vp(u_2))$$and rewrite \eqref{e3.13bb} as
\begin{equation}\label{e3.17ae}
u={\rm div}\ F_\vp-\vp\lbb(\beta(u_1)-\beta(u_2))+f.\end{equation}
  By \eqref{e3.12c},   it follows that $F_\vp\in L^2(\rr^d)$ and by \eqref{e3.17ae} that ${\rm div}\,F_\vp\in L^2(\rrd)$.   
  We set $\Lambda_\delta=\calx_\delta(\beta(u_1)-\beta(u_2))$. 
  Since $\Lambda_\delta\in H^1(\rr^d)$, it follows that $\Lambda_\delta{\rm div}\,F_\vp\in L^1$ and so, by \eqref{e3.17ae}, we~have  
  \begin{equation}\label{e3.17e}
  \barr{ll}
  \dd\int_{\rr^d}u\Lambda_\delta dx\!\!\!
  &=-\dd\int_{\rr^d}F_\vp\cdot\nabla
  \Lambda_\delta dx\vsp
  &-\,\vp\lbb\dd\int_{\rr^d}(\beta (u_1)-\beta (u_2))\Lambda_\delta dx
   +\dd\int_{\rr^d}f\Lambda_\delta dx\vsp
  &=-\dd\int_{\rr^d}(F_\vp\cdot\nabla (\beta(u_1)-\beta(u_2))\calx'_\delta(\beta(u_1)-\beta(u_2))dx\vsp
  &-\vp\lbb\dd\int_{\rr^d}(\beta(u_1)-\beta(u_2))\calx_\delta(\beta(u_1)-\beta(u_2))dx
   +\dd\int_{\rr^d}f\Lambda_\delta dx.\earr\hspace*{-10mm}
  \end{equation}
 We set
\begin{equation}\label{e3.20f}
\hspace*{-5mm}\barr{ll}
I^1_{\delta}\!\!\!&=
\dd\int_{\rr^d}E_\vp (b^*_\vp(u_1)-b^*_\vp(u_2))
\cdot\nabla\Lambda_\delta dx\vsp
&=\dd\int_{\rr^d}E_\vp (b^*_\vp(u_1)-b^*_\vp(u_2))
\cdot\nabla (\beta(u_1)-\beta(u_2)) \calx'_\delta(\beta(u_1)-\beta(u_2))dx\vsp
&=\dd\frac1\delta
\int_{[|\beta(u_1)-\beta(u_2)|\le\delta]}E_\vp (b^*_\vp(u_1)-b^*_\vp(u_2))\cdot\nabla (\beta(u_1)-\beta(u_2))dx.\earr\end{equation}
Since $|E_\vp|\in L^\9\cap L^2$  and, by Hypothesis (i),
$$|b^*_\vp(u_1)-b^*_\vp(u_2)|\le {\rm Lip}(b^*_\vp)|u_1-u_2|
\le\frac1\gamma\ {\rm Lip}(b^*_\vp)|\beta(u_1)-\beta(u_2)|,$$
   it follows that
$$\barr{l}
	\dd\lim_{\delta\to0} \frac1\delta\int_{[|(\beta(u_1)-\beta(u_2))|\le\delta]}|E_\vp
(b^*_\vp(u_1)-b^*_\vp(u_2))\cdot\nabla (\beta(u_1)-\beta(u_2))|dx\\
\dd\le\frac1\gamma\ {\rm Lip}(b^*_\vp)|E_\vp |_2
\lim_{\delta\to0}\(\int_{[|\beta(u_1)-\beta(u_2)|\le\delta]} |\nabla(\beta(u_1)-\beta(u_2))|^2dx\)^{\frac12} \!\!=0.\earr$$
 This yields
\begin{equation}\label{e3.13bbbbb}
 \lim_{\delta\to0}I^1_{\delta}=0,
\end{equation}
because $\nabla (\beta(u_1)-\beta(u_2))(x)=0$, a.e. on $[x\in\rr^d;\beta(u_1(x))-\beta(u_2(x)){=}0]$.
On the other hand, since $\calx'_\delta\ge0$, we have
\begin{equation}\label{e3.23c}
\dd\int_{\rr^d}\nabla(\beta(u_1)-\beta(u_2))
\cdot\nabla (\beta(u_1)-\beta(u_2))\calx'_\delta
(\beta(u_1)-\beta(u_2))\,dx\ge0.\end{equation}
  By \eqref{e3.17e}-\eqref{e3.23c},  since $|\Lambda_\delta|\le1,$ we get
$$\lim_{\delta\to0}\int_{\rr^d} u\calx_\delta(\beta(u_1)-\beta(u_2))dx\le \int_{\rr^d} |f|\,dx$$and, since $u\calx_\delta(\beta(u_1)-\beta(u_2))\ge0$ and $\calx_\delta\to{\rm sign}$ as $\delta\to0$, by Fatou's lemma  this yields
\begin{equation}\label{e3.18c}
|u|_1\le|f|_1,\end{equation}as claimed.

Next, for $f$ arbitrary in $L^1$,   consider a sequence $\{f_n\}\subset L^2$ such that $f_n\to f$ strongly in $L^1.$
Let $\{u^n_\vp\}\subset L^1\cap L^2$ be the corresponding solutions to \eqref{e3.8} for $0<\lbb<\lbb_\vp$. We have, for all $m,n\in\mathbb{N}$,
\begin{equation}\label{e3.18cc}
u^n_\vp-u^m_\vp+\lbb((A_0)_\vp u^n_\vp-(A_0)_\vp u^m_\vp)=f_n-f_m.\end{equation}Taking into account  \eqref{e3.18c}, we obtain by the above equation that
$$|u^n_\vp-u^m_\vp|_1\le|f_n-f_m|_1,\ \ff n,m\in\nn.$$Hence, for $n\to\9$, we have
 $u^n_\vp\to u_\vp(\lbb,f)\mbox{ in }L^1.$
Now, \eqref{e3.18cc} implies that $(A_0)_\vp u^n_\vp\to v\mbox{ in }L^1.$ Since $((A_0)_\vp,D((A_0)_\vp))$ is closed on $L^1$, we conclude that $u_\vp(\lbb,f)\in D((A_0)_\vp)$ and that
\begin{equation}\label{e3.25v}
u_\vp(\lbb,f)+\lbb (A_0)_\vp u_\vp (\lbb,f)=f,
\end{equation}which proves \eqref{e3.6a}    for $\lbb<\lbb_\vp$. Moreover, by \eqref{e3.18c}, we have
\begin{equation}\label{e3.13}
|u_\vp(\lbb,f_1)-u_\vp(\lbb,f_2)|_1\le|f_1-f_2|_1,\ \ff f_1,f_2\in L^1.\end{equation} 
By Proposition 3.3 in \cite{1}, p. 99, it follows that  $R(1+\lbb (A_0)_\vp)=L^1,$ $ \ff\lbb>0,$ and, therefore, \eqref{e3.25v} holds for all $\lbb>0$ if $f\in L^1$.
We set $J^\vp_\lbb(f)=u_\vp(\lbb,f).$ Then, by  \eqref{e3.25v}, \eqref{e3.13}, it follows that \eqref{e3.6a}, \eqref{e3.6aaa}, \eqref{e3.6aaaa} are satisfied. Since $u_\vp=J^\vp_\lbb(f)$ is for $f\in L^1\cap L^2$ the solution to \eqref{e3.7}, it follows \eqref{e3.6aa} for all $f\in L^1\cap L^2$ and so by density for all $f\in L^1$. We also note that, by \eqref{e3.7},
\begin{equation}
 \barr{r}
 \dd\int_{\rr^d} J^\vp_\lbb(f)dx=\int_{\rr^d} f\,dx
-\vp\lbb\int_{\rr^d} \beta(J^\vp_\lbb(f))dx,\\ \ff f\in L^1\cap L^2,\ \lbb>0,\earr\label{e3.14}\end{equation}and so \eqref{e3.6aaaaaa} follows for all $f\in L^1\cap L^2$ and so, by \eqref{e3.13} for all $f\in L^1$. Note also that there exists $\wt\lbb_1$ independent of $\vp$ such that, for all $\lbb\in(0,\wt\lbb_1)$ and $f\in L^1\cap L^2$,
\begin{equation}
J^\vp_\lbb(f)\ge0,\ \mbox{a.e. in }\rr^d\mbox{ if }f\ge0,\mbox{ a.e. in }\rr^d. \label{e3.15}
\end{equation}(The latter follows by multiplying \eqref{e3.7}, where $u=u_\vp$, with sign $u^-_\vp$ and integrating over $\rr^d$.)

Next, we show \eqref{e3.3}. Fix $\lbb<\lbb_0=\min(\lbb_1,\wt\lbb_1)$, with $\lbb_1$ as in \eqref{e3.13a},  and let $f\in L^1\cap L^2$. If $u_\vp=u_\vp(\lbb,f)$, by \eqref{e3.13a}, it follows that $\{u_\vp\}$ is bounded in $H^1(\rr^d)$ and $\{\beta(u_\vp)\}$ is bounded in $H^1(\rr^d)$.  Clearly,  $u_\vp(f)=0$ if $f\equiv0$, hence \eqref{e3.13}  implies that $\{u_\vp\}$ is bounded in $L^1$.  Hence, along a subsequence, again denoted $\{\vp\}\to0$, we~have\newpage
\begin{equation}\label{e3.18}
\barr{rcll}
u_\vp&\!\!\!\longrightarrow\!\!\!&u&\mbox{weakly in $H^1(\rr^d)$, strongly in $L^2_{\rm loc}(\rr^d)$,}\vsp
\beta(u_\vp)&\!\!\!\longrightarrow\!\!\!&\beta(u)
&\mbox{weakly in $H^1(\rr^d)$ and strongly in $L^2_{\rm loc}(\rr^d)$},\vsp	 \D\beta(u_\vp)&\!\!\!\longrightarrow\!\!\!&\D\beta(u)&\mbox{weakly in $H^{-1}(\rr^d)$},\earr\end{equation}
and, by Hypothesis (ii) and \eqref{e3.2a},
\begin{equation}\label{e3.30prim}
\barr{rcll}
 b^*_\vp(u_\vp)&\longrightarrow&b(u)u
&\mbox{strongly in $L^2_{\rm loc}(\rr^d)$.}\earr \end{equation}This yields
\begin{equation}\label{e331}
E_\vp b^*_\vp(u_\vp)\to Eb(u)u\mbox{ strongly in }L^2_{\rm loc}(\rr^d).\end{equation}
Passing to the limit in \eqref{e3.7}, we obtain
\begin{equation}\label{e3.32i}
u-\lbb\Delta\beta(u)+\lbb\ {\rm div}(Eb(u)u)=f\mbox{ in }\cald'(\rr^d),
\end{equation}
 where $u=u(\lbb,f)\in H^1(\rr^d)$.
By \eqref{e3.13} and \eqref{e3.18}, it follows via Fatou's lemma that
\begin{equation}\label{e3.21}
|u(\lbb,f_1)-u(\lbb,f_2)|_1\le|f_2-f_2|_1,\ \ff f_1,f_2\in L^2\cap L^1,
\end{equation}and hence (since $u(\lbb,f)=0$ if $f\equiv0$) $u_1(\lbb,f),u_2(\lbb,f)\in L^1\cap L^2$, if $f\in L^1\cap L^2.$
 In particular, $u(\lbb,f)\in D(A_0)$ and
	\begin{equation}\label{e3.34b}
	u(\lbb,f)+\lbb A_0u(\lbb,f)=f,\ \ff f\in L^1\cap L^2.\end{equation}Now, let $f\in L^1$ and $f_n\in L^1\cap L^2$, $n\in\nn$, such that $f_n\to f$ in $L^1$. Then, by \eqref{e3.21}, $u(\lbb,f_n)\to u=u(\lbb,f)$ in $L^1$ and, therefore, since each $u(\lbb,f_n)$ satisfies \eqref{e3.34b}, we conclude that $u(\lbb,f)\in D(A_0)$ and that $u$ also satisfies \eqref{e3.34b}, and so \eqref{e2.2} follows for all $\lbb\in(0,\lbb_0)$. Again by Proposition 3.3 in \cite{1}, p.~99, \eqref{e2.2} and \eqref{e3.21} extend to all $\lbb>0$.
	
	We define $J_\lbb:L^1\to L^1$ as $J_\lbb(f)=u(\lbb,f)$ and, by \eqref{e3.21},  \eqref{e2.3} follows. Moreover, letting $\vp\to0$ in \eqref{e3.6aa}--\eqref{e3.6aaaaa}, it follows that $J_\lbb$ satisfies \eqref{e2.3a}--\eqref{e2.3} and \eqref{e2.4}, \eqref{e2.5}, as claimed. 
	
Clearly, by \eqref{e3.18},
\begin{equation}\label{e3.36x}
u_\vp\to u=u(\lbb,f)=J_\lbb(f)\mbox{ in }L^1_{\rm loc},
\end{equation}	
	for $0<\lambda<\lambda_0$. (Here, $u_\vp=J^\vp_\lbb(f)=(I+\lbb A_\vp)\1f$.)
	
To prove that \eqref{e3.3}, that is that \eqref{e3.36x}   holds in $L^1$,   we shall prove first the following lemma, which has an intrinsic interest and where we use Hypothesis (iv) for the first time.

\begin{lemma}\label{l3.2} Assume that Hypotheses {\rm(i)-(iv)}   hold, and let  $u_0\in  \calm\cap L^2.$
	\begin{itemize}
		\item[\rm(a)] We have
		\begin{equation}
		\label{e3.39'}
		\sup_{\vp\in(0,1)}\int_\rrd|u_\vp|\Phi\,dx<\9.
		\end{equation}
	 \item[\rm(b)] Assume that ${\rm div}\,E\in L^\9$. Then, for all $\lbb\in(0,\lbb_0)$,
\begin{equation}\label{e3.37x}
\|(I+\lbb  A_\vp)\1u_0\|\le
  \|u_0\|+\rho_\vp\lbb |u_0|_1,\end{equation}
  where  $\rho_\vp=\gamma_1(m+1)|\Delta\Phi|_\9+\gamma_1
m(m+3)\vp|E|^2_\9+|b|_\9(1+m)^2|E|^2_\9.$\end{itemize}
\end{lemma}
\n{\bf Proof.} By approximation also in (b), we may restrict to the case   $u_0\!\in\!\calm\cap L^2.$\break   If we  multiply equation \eqref{e3.25v} by $\vf_\nu\calx_\delta(\beta(u_\vp))$, where $u_\vp=(I+\lbb (A_0)_\vp)\1u_0=(I+\lbb A_\vp)^{-1}u_0$,    $\vf_\nu(x)=\Phi_\vp(x)\exp(-\nu\Phi_\vp(x))$ and integrate over $\rr^d$, we get, since $\calx'_\delta\ge0,$
\begin{equation}\label{e3.36av}
\barr{l}
 \dd\int_{\rr^d}u_\vp\calx_\delta
 (\beta(u_\vp))\vf_\nu\,dx
 \le-\lbb \dd\int_{\rr^d}\nabla\beta(u_\vp)
 \cdot\nabla(\calx_\delta(\beta(u_\vp))\vf_\nu)dx\\
  \qquad+\lbb \dd\int_{\rr^d}E_\vp b^*_\vp(u_\vp)\cdot
 \nabla(\calx_\delta(\beta(u_\vp))\vf_\nu)dx
 +\dd\int_{\rr^d}|u_0|\vf_\nu dx
\\
 \qquad\le-\lbb\dd\int_{\rr^d}\nabla\beta(u_\vp)
 \cdot\nabla\vf_\nu\calx_\delta(\beta(u_\vp))dx\\
 \qquad
 +\lbb \dd\int_{\rr^d}E_\vp b^*_\vp(u_\vp)\cdot\nabla
 \beta(u_\vp)\calx'_\delta(\beta(u_\vp))\vf_\nu dx\vsp
 \qquad+\lbb \dd\int_{\rr^d}(E_\vp\cdot\nabla\vf_\nu)b^*_\vp(u_\vp)\calx_\delta(\beta(u_\vp))dx
 +\dd\int_{\rr^d}|u_0|\vf_\nu dx.\earr\end{equation}
Letting $\delta\to0$, we get as above
\begin{equation}\label{e3.37za}
\barr{l}
 \dd\int_{\rr^d}|u_\vp|\vf_\nu dx
 \le-\lbb\dd\int_{\rr^d}\nabla|\beta(u_\vp)|\cdot\nabla\vf_\nu dx\\
 \quad
 +\overline{\lim\limits_{\delta\to0}}\dd\frac\lbb\delta
 \dd\int_{[|\beta(u_\vp)|\le\delta]}
 |E_\vp|\,|b^*_\vp(u_\vp)|\,|\nabla\beta(u_\vp)|\vf_\nu dx\\
 \quad
 +\lbb\dd\int_{\rr^d}{\rm sign}\,u_\vp b^*_\vp(u_\vp)
 E_\vp\cdot\nabla\vf_\nu\,dx+
 \dd\int_{\rr^d}\!|u_0|\vf_\nu dx\vsp\quad
 \le\lbb\dd\int_{\rr^d}
 (|\beta(u_\vp)|
 \Delta\vf_\nu
 +|b^*_\vp(u_\vp)|
  \,|\nabla\Phi_\vp\cdot\nabla\vf_\nu|)dx
\! +\!\!\dd\int_{\rr^d}\!\!|u_0|\vf_\nu dx,
 \earr\end{equation}because $|b^*(u_\vp)|\le C|u_\vp|\le \frac C\gamma\,|\beta(u_\vp)|$, a.e. in $\rr^d$, and so
$$\frac1\delta\int_{[|\beta(u_\vp)|\le\delta]}|E_\vp|\,|b^*(u_\vp)|\,|\nabla\beta(u_\vp)|\vf_\nu dx\le \frac C\gamma\,|E_\vp|_2\(\int_{[|\beta(u_\vp)|\le\delta]}|\nabla\beta(u_\vp)|^2dx\)^{\frac12}$$and
$$\lim_{\delta\to0}\int_{[|v|\le\delta]}|\nabla v|^2dx=0,\ \ \ff v\in H^1(\rr^d).$$
We have
\begin{eqnarray}
 \nabla\vf_\nu(x)&=&(1-\nu\Phi_\vp)  \nabla\Phi_\vp\exp(-\nu\Phi_\vp),\label{e3.37a}\\
 \Delta\vf_\nu(x)&=&((1-\nu\Phi_\vp)
 \Delta\Phi_\vp-2\nu|\nabla\Phi_\vp|^2
 +\nu^2\Phi_\vp|\nabla\Phi_\vp|^2)\exp(-\nu\Phi_\vp),\qquad\label{e3.37aa}\\
 \Delta\Phi_\vp&=&-\,{\rm div}\ E_\vp
 =(1-m\vp\Phi(1+\vp\Phi)^{-1})(1+\vp\Phi)^{-m}\Delta\Phi\label{e340a}\\
 &&+\,m\vp((m+1)\vp\Phi
 (1+\vp\Phi)^{-1}-2) (1+\vp\Phi)^{-(m+1)}|E|^2.\nonumber
 \end{eqnarray}
 Then, letting $\nu\to0$, since $\beta(u_\vp),$ $\vp\in(0,1)$, is bounded in $L^1\cap L^2$, we get by \eqref{e3.37za}, \eqref{e340a} and Hypothesis (iii) that
 $$\sup_{\vp\in(0,1)}\int_\rrd|u_\vp|\Phi\,dx<\9,$$and assertion (a) follows.   If ${\rm div}\,D\in L^\9$, we additionally get from \eqref{e3.37za}~that
  	$$\|u_\vp\|
  	 \le\|u_0\|+ \lbb \gamma_1|\Delta\Phi_\vp|_\9 |u_0|_1+\lbb|b|_\9|u_0|_1|\nabla\Phi_\vp|^2_2,\ \ff\vp>0.$$By \eqref{e340a}, we have
  	 \begin{equation}\label{e3.45prim}
  	 |\Delta\Phi_\vp(x)|\le(m+1)|\Delta\Phi(x)|+m(m+3)\vp|E|^2(x)\mbox{ for a.e. }x\in\rrd,\end{equation}
  	and this, together with \eqref{e3.4x},  yields
  	\eqref{e3.37x}, as claimed.

\begin{remark}\label{r3.3}\rm If, as in Remark \ref{r2.4}, we replace (iv), $\calm$, $\|\cdot\|$ and $\rho$ by (iv)$'$ (see Remark \ref{r2.4}), $\calm_2$, $\|\cdot\|_2$ and $\wt\rho$, respectively, we can prove a complete analogue of Lemma \ref{l3.2} by the same arguments. One only has to replace $\vf_\nu$ by the function $\wt\vf_\nu(x)=|x|^2e^{-\nu|x|^2}$ in the above proof. Once one has this analogue of Lemma \ref{l3.2}, the proofs below can easily be adjusted to this case.
\end{remark}
 \n{\bf Proof of \eqref{e3.3}.} By \eqref{e3.39'} and Hypothesis (iv),  it follows that, if $f\in\calm\cap L^2$, then we have, for all $\lbb\in(0,\lbb_0)$ and $\vp\in(0,1)$, $N>0$,
 $$\int_{\{\Phi\ge N\}}|(I+\lbb A_\vp)\1 f|dx\le \frac 1N\|(I+\lbb A_\vp)\1f\|
 \le\frac CN.$$
 Recalling \eqref{e3.36x} and that $\{\Phi\le N\}$ is compact,  the latter implies that, if $f\in\calm\cap L^2,$ then $\lim_{\vp\to0}|u_\vp-u|_1=0,$ i.e.,
\begin{equation}\label{e3.22a}
\lim_{\vp\to0}(I+\lbb A_\vp)^{-1}f=(I+\lbb A)^{-1}f\mbox{ in }L^1,\ \ff f\in \calm\cap L^2.
\end{equation}
Since $L^2\cap\calm$ is dense in $L^1$ and $(I+\lbb A_\vp)\1,$ $\vp>0$, are equicontinuous,   \eqref{e3.3} follows.\bigskip
 
\noindent{\bf Proof of \eqref{e2.3d}.}
Let $f\in C^\9_0(\rrd)$ and $u_\lbb=J_\lbb(f)\in D(A),$ $\lbb>0$. Since $D(A)\subset D(A_0)$, we have
	\begin{equation}\label{e3.49}
	u_\lbb+\lbb A_0 u_\lbb=f,\end{equation}where $u_\lbb\in L^1\cap L^\9$, $|u_\lbb|_1\le|f|_1,$ and, by Lemma 3.1 in \cite{6a},
	\begin{equation}
	\label{e3.50}
 \sup_{\lbb\in(0,\lbb_0)}|u_\lbb|_\9=C_\9<\9.
	\end{equation} By \eqref{e3.13a}, we also have
	\begin{equation}
	\label{e3.51}
	\sup_{\lbb\in(0,\lbb_0)}|u_\lbb|^2_2=C_2<\9,
	\end{equation}for some $\lbb_0>0.$ Taking into account \eqref{e3.50} and that $b^*(r)\equiv b(r)r$ is locally Lipschitz, it follows as in the proof of Lemma \ref{l3.2} that 
		\begin{equation}
	\label{e3.52}
	\sup_{\lbb\in(0,\lbb_0)}\int_\rrd|u_\lbb(x)|\Phi(x)dx<\9.
	\end{equation}Next, by \eqref{e3.49}, we see that since $A_0u_\lbb\in L^2$, we have
	$$\<A_0u_\lbb,u_\lbb\>_2+\lbb|A_0u_\lbb|^2_2=\<A_0u_\lbb,f\>_2\le|A_0u_\lbb|_2|f|_2.$$This yields
	$$\<\nabla\beta(u_\lbb),\nabla u_\lbb\>_2\le
	\<E,b^*(u_\lbb)\nabla u_\lbb\>_2+
	\<\nabla\beta(u_\lbb),\nabla f\>_2-(E,b^*(u_\lbb)\nabla f)_2$$and so, by Hypotheses (i)--(ii) we get, for $\delta>0$,
	$$\gamma|\nabla u_\lbb|^2_2\le\delta(1+\gamma^2_1)|\nabla u_\lbb|^2+\dd\frac1\delta(|E|^2_\9|b|^2_\9|u_\lbb|^2_2+|\nabla f|^2_2)
	+|E|_\9|b|_\9|u_\lbb|_2|\nabla f|_2.$$
This yields
		\begin{equation}
	\label{e3.53}
	|\nabla u_\lbb|^2_2\le K_\delta(\gamma-\delta(1+\gamma^2_1))\1.
	\end{equation}By \eqref{e3.49}--\eqref{e3.53}, it follows
	$$\lbb A_0u_\lbb\to0\mbox{\ \ in }H\1\mbox{ as }\lbb\to0$$and, therefore, $u_\lbb\to f$ in $H\1$ as $\lbb\to0$ and so, by \eqref{e3.53}, we have on a subsequence $\{\lbb\}\to0$
	$$u_\lbb\to f\mbox{\ \ in }L^2_{\rm loc}\subset L^1_{\rm loc}.$$\newpage \n Then, by \eqref{e3.52}, we infer that for $\lbb\to0$, $u_\lbb\to f$ in $L^1$ and so $f\in\ov{D(A)}$. Hence, $C^\9_0(\rrd)\subset\ov{D(A)}$ and so \eqref{e2.3d} follows.
	
We note that, similarly, it follows that
		\begin{equation}
	\label{e3.54}
	\ov{D(A_\vp)}=L^1.
	\end{equation}
 This completes the proof of Proposition \ref{p2.1}.

\bk\n{\bf Proof of Proposition \ref{p2.3}.} By Lemma \ref{l3.1} and \eqref{e3.37x} in Lemma \ref{l3.2}, we have, for $\lbb\in(0,\lbb_0)$,   and $\delta>0$,
$$\|(I+\lbb A)^{-1}u_0\|\le
   \|u_0\|+\rho  \lbb|u_0|_1,\ \ff u_0\in\calm  .$$
This yields
$$ \|(I+\lbb A)^{-n}u_0\|\le
\|u_0\|+n\lbb\rho|u_0|_1,\  \ff n\in\nn,$$
and so, by \eqref{e2.10}, we get
\begin{equation}\label{e3.43}
\|S(t)u_0\|\le
 \|u_0\|+\rho t|u_0|_1,\ \ff t\ge0,\ u_0\in\calm,\end{equation}
as claimed.
\section{The $H$-theorem}\label{s4}
\setcounter{equation}{0}

Let $S(t)$ be the continuous semigroup of contractions defined by \eqref{e2.12}.
A  lower semicontinuous function $V:L^1\to(-\9,\9]$ is said to be a {\it Lyapunov function} for $S(t)$ (equivalently, for equations \eqref{e1.1} or \eqref{e2.6}) if
$$V(S(t)u_0)\le V(S(s)u_0),\mbox{ for }0\le s\le t<\9,\ u_0\in L^1.$$
(See, e.g.,   \cite{9}.)

In the following, we shall restrict the semigroup to the probability density set $\calp$ (see \eqref{e2.16}).
For each $u_0\in\calp$, consider the $\oo$-limit set
$$\oo(u_0)=\{w=\lim S(t_n)u_0\mbox{ in }L^1_{\rm loc}\mbox{ for some }\{t_n\}\to\9\}.$$
Our aim here is to construct a Lyapunov function for $S(t)$,   to prove that $\oo(u_0)\ne\emptyset$ and also that every $u_\9\in\oo(u_0)$ is an equilibrium state of equation \eqref{e1.1}, that is, $Au_\9=0$.
To this end, we shall assume that, besides (i)-(iv),   Hypothesis (v) also holds.

Consider the function $\eta\in C(\rr)$,
\begin{equation}\label{e4.2}
\eta(r)=-\int^r_0d\tau\int^1_\tau\frac{\beta'(s)}{sb(s)}\,ds,\ \ff r\ge0,\end{equation}
and define the function \mbox{$V:D(V)=\calm_+=\{u\in \calm; u\ge0,\mbox{ a.e. on }\rr^d\}\to\rr$}   \begin{equation}\label{e4.3}
V(u)=\int_{\rr^d}\eta(u(x))dx+\int_{\rr^d}\Phi(x)u(x)dx=-\wt S[u]+F[u]. \end{equation}
Since, by (i), (iv) and (v),
\begin{equation}\label{e4.1av}
\frac \gamma{r|b|_\9}\le\frac{\beta'(r)}{rb(r)}\le\frac{\gamma_1}{rb_0},\ \ff r>0,
\end{equation}
we have
\begin{equation}\label{e4.3e}
\barr{r}
\dd\frac{\gamma_1}{b_0}\,\one_{[0,1]}(r)r(\log r-1)
+\frac{\gamma}{|b|_\9}\,\one_{(1,\9)}(r)
r(\log r-1)\le\eta(r)\vsp
\le\dd\frac{\gamma}{|b|_\9}\,\one_{[0,1]}(r)r(\log r-1)+
\dd\frac{\gamma_1}{b_0}\,\one_{(1,\9)}(r)r(\log r-1).\earr\end{equation}
We also have   that $\eta\in C([0,\9)),$ $\eta\in C^2((0,\9)),$ $ \eta''\ge0.$
 Since $\Phi$ is Lipschitz, hence of at most linear growth, $F[u]$ is well-defined and finite  if $u\in\calm$. Furthermore, exactly as in \cite{11b}, p.~16, one proves   that $(u\log u)^-\in L^1$ if $u\in D(V)$. Hence $\wt S[u]$ is well-defined and $-\wt S[u]\in(-\9,\9]$   because of  \eqref{e4.3e} and thus $V(u)\in(-\9,\9]$ for all $u\in D(V)$. We define $V=\9$ on $L^1\setminus D(V)$. Then, obviously, $V:L^1\to(-\9,\9]$ is convex and $L^1_{\rm loc}$-lower semicontinuous on balls in $\calm$, as easily follows by \eqref{e4.3e} from \eqref{e4.5z} below.  If, in addition, $(u\log u)^+\in L^1$, then, again by \eqref{e4.3e},  we have that $\wt S[u]\in(-\9,\9)$ and also $V$ is real-valued.
The function (see \eqref{e1.8})
$$\wt S[u]=-\int_{\rr^d}\eta(u(x))dx,\ u\in\calp,$$is called in the literature (see, e.g., \cite{7}, \cite{10}) the   entropy of the system, while $F[u]$ is the mean field energy. In fact, according to the general theory of thermostatics (see \cite{8}), the functional $\wt S=\wt S[u]$ is a generalized entropy because its kernel $-\eta$ is a strictly  concave continuous functions on $(0,\9)$ and $\lim\limits_{r\downarrow0}\eta'(r)=+\9.$ In the special case $\beta(s)\equiv s$ and $b(s)\equiv1$, $\eta(r)\equiv r(\log r-1)$ and so $\wt S[u]-1$ reduces to the classical Boltzmann-Gibbs entropy.

As in \cite{11b} (formula (15)), one proves that, for $\alpha\in\left[\frac m{m+1},1\right)$, where $m$ is as in assumption (iv),
\begin{equation}\label{e4.5z}
\int_{\{\Phi\ge R\}}|\min(u\log u,0)|dx\le C_\alpha\(\int_{\{\Phi\ge R\}}\Phi^{-m}dx\)^{1-\alpha}\|u\|^\alpha,\end{equation}for all $R>0$.
Indeed, obviously, for every $\alpha\in(0,1)$, there exists $C_\alpha\in(0,\9)$ such that
  $(r\log r)^-\le C_\alpha r^\alpha\mbox{\ \ for }r\in[0,\9).$
Hence, the left hand side of \eqref{e4.5z} by H\"older's inequality is dominated by
$$C_\alpha\left(\int_{\{\Phi\ge R\}}u\Phi dx\)^\alpha\(\int_{\{\Phi\ge R\}}
\Phi^{-\frac\alpha{1-\alpha}}dx\)^{1-\alpha}.$$Therefore, for $\alpha\in\left[\frac m{m+1},1\right)$, we obtain \eqref{e4.5z} since $\Phi\ge1$ and \eqref{e4.5z} yields
\begin{equation}\label{e4.3am}
V(u)\ge-C(\|u\|+1)^\alpha,\ \ff u\in D(V).
\end{equation}
We also consider   the function $\Psi:D(\Psi)\subset L^1\to[0,\9)$ defined by
\begin{eqnarray}
\Psi(u)&=&\int_{\rr^d}\left|\frac{\beta'(u)\nabla u}{\sqrt{ub(u)}}-E\sqrt{ub(u)}\right|^2 dx,\label{e4.4}\\[2mm]
D(\Psi)&=&\{u\in L^1\cap W^{1,1}_{\rm loc}(\rr^d);\ u\ge0,\ \Psi(u)<\9\}.\label{e4.5}
\end{eqnarray}
We extend $\Psi$ to all of $L^1$ by $\Psi(u)=\9$ if $u\in L^1\setminus D(\Psi)$.
Since $\nabla u=0$, a.e. on $\{u=0\}$, we set here and below
$$\frac{\nabla u}{\sqrt{u}}=0\mbox{\ \ on }\{u=0\}.$$
Theorem \ref{t4.1} is the main result and, as mentioned earlier, can be viewed as the $H$-theorem for NFPE \eqref{e1.1}.

\begin{theorem}\label{t4.1} Assume that Hypotheses {\rm(i)-(v)} and \eqref{e2.21a} hold. Then the function $V$ defined by \eqref{e4.3} is a Lyapunov function for $S(t)$, that is, for $D_0(V)=D(V)\cap \{V<\9\}$ $(=\{u\in D(V);\,u\log u\in L^1\})$,
\begin{equation}\label{e4.8prim}
\barr{c}
S(t)u_0\in D_0(V),\,\ff t\ge0,u_0\in D_0(V)\mbox{ and }\vsp V(S(t)u_0)\le V(S(s)u_0),\,\ff u_0\in D_0(V),  0\le s\le t<\9.\earr\end{equation}  Moreover, we have, for all $u_0\in D_0(V)$,
\begin{equation}
 V(S(t)u_0)+\dd\int^t_s\Psi(S(\sigma)u_0)d\sigma\le V(S(s)u_0)\mbox{ for $0\le s\le t<\9$}.\label{e4.6}\end{equation}
 In particular, $S(\sigma)u_0\in D(\Psi)$ for  a.e. $\sigma\ge0.$ Furthermore, there exists $u_\9\in\oo(u_0)$ $($see \eqref{e1.4}$)$ such that $u_\9\in D(\Psi)$, $\Psi(u_\9)=0$ and, for any such~a $u_\9$, we have either $u_\9=0$ or $u_\9>0$ a.e. In the latter case,
 \begin{equation}\label{e4.11i}
 u_\9=g^{-1}(-\Phi+\mu)\mbox{ for some }\mu\in\rr,\end{equation}
\begin{equation}\label{e4.12i}
 g(r)=\int^r_1\frac{\beta'(s)}{sb(s)}\ ds,\ r>0.
 \end{equation}
\end{theorem}

Moreover, by \eqref{e4.3}, \eqref{e4.6}, we see that the entropy of the semiflow $u(t)=S(t)u_0$  is evolving according to the law
$$\wt S[u(t)]\ge \wt S[u(s)]+\dd\int_{\rr^d}\Phi(x)(u(t,x)-u(s,x))ds
+\int^t_s\Psi(u(\sigma))d\sigma,$$for all $ 0\le s\le t<\9. $

\begin{remark}\label{r4.2}\rm We note that \eqref{e2.21a} holds if ${\rm div}\,D\in L^\9$ (see Propopsition~\ref{p2.3}) or  if Hypothesis (vi) holds (see Lemma \ref{l6.2}  below).\end{remark}

\section{Proof of Theorem \ref{t4.1}}\label{s5}
\setcounter{equation}{0}

In the following, we approximate $V:L^1\to(-\9,\9]$ by the functional $V_\vp$ defined by
$$\barr{c}
V_\vp(u)=\dd\int_{\rr^d}(\eta_\vp(u(x))+\Phi_\vp(x)u(x))dx,\ \ff u\in D(V),\vsp
V_\vp(u)=\9\mbox{\ \ if }u\in L^1\setminus D(V),\earr$$where
$ \eta_\vp(r)=-\int^r_0d\tau\int^1_\tau\,
\frac{\beta'(s)}{b^*_\vp(s)+\vp^{2m}}\,ds,\ r\ge0,\ \vp>0.$
Clearly, $\eta_\vp\to\eta$ as $\vp\to0$ locally uniformly.
We note that $V_\vp$ is convex, and $L^1_{\rm loc}$-lower semicontinuous on every ball in $\calm$. Furthermore, there exists $C>0$ such that, for all $\vp\in(0,1]$, we have  $|\eta_\vp(u)|\le C(1+|u|^2).$ This implies that $V_\vp<\9$ on $L^2$ and $V_\vp(u)\to V(u)$ as $\vp\to0$ for all $u\in D(V)\cap L^2$ and by the generalized Fatou  lemma that $V_\vp$ is lower semicontinuous on $L^2$.  We set 
$$V'_\vp(u)=\eta'_\vp(u)+\Phi_\vp,\ \ff u\in D(V)\cap L^2.$$\newpage\n It is easy to check that $V'_\vp(u)\in\partial V_\vp(u)$ for all $u\in D(V)\cap L^2$, where $\partial V_\vp$ is the subdifferential of $V_\vp$ on $L^2$.
As regards the function $\Psi$ defined by \eqref{e4.4}-\eqref{e4.5}, we have
\begin{lemma}\label{l5.1} We have
\begin{equation}\label{e5.3}
D(\Psi)=\{u\in L^1;\ u\ge0,\ \sqrt{u}\in W^{1,2}(\rr^d)\},\end{equation}
\begin{equation}\label{e5.4}
\|\sqrt{u}\|_{W^{1,2}(\rr^d)}\le C(\Psi(u)+1),\ \ff u\in D(\Psi),\end{equation} where $C\in (0,\9)$ is independent of $u$. Furthermore, $\Psi$ is  $L^1_{\rm loc}$-lower semi\-continuous on $L^1$-balls.
\end{lemma}
\noindent{\bf Proof.} By \eqref{e4.4}, taking into account   (i), (ii), we have
\begin{equation}\label{e52a}
 \barr{ll}
\dd\gamma|b|^{-1}_\9 \int_{\rr^d}
\frac{|\nabla u|^2}u\,dx
&\le \dd\int_{\rr^d}
\frac{|\beta'(u)|^2\cdot|\nabla u|^2}{ub(u)}\,dx\vsp
&\le2\Psi(u)+2\!\!\dd\int_{\rr^d}|E|^2ub(u)dx<\9,\  \ff u\in D(\Psi).\earr\hspace*{-5mm}\end{equation}
This yields \eqref{e5.3} and \eqref{e5.4} since $\nabla(\sqrt{u})=\frac12\ \frac{\nabla u}{\sqrt{u}}$ and (v) holds.
To show the lower semicontinuity of $\Psi$, we rewrite it as
\begin{equation}\label{e5.5}
\Psi(u)=\int_{\rr^d}|\nabla j(u)-E\sqrt{ub(u)}|^2dx,\ u\in D(\Psi),\end{equation}where
\begin{equation}\label{e5.4c}
j(r)=\int^r_0\frac{\beta'(s)}{\sqrt{sb(s)}}\ ds,\ r\ge0.\end{equation}
Clearly,
\begin{equation}\label{e55a}
0\le j(r)\le\frac{2\gamma_1}{\sqrt{b_0}}\ \sqrt{r}.\end{equation}
Let $\{u_n\}\subset L^1$ and $\nu>0$ be such that $\sup\limits_n|u_n|_1<\9$ and
\begin{equation}\label{e5.6}
\Psi(u_n)\le\nu<\9,\ \ff n,\end{equation}
\begin{equation}\label{e5.7}
u_n\longrightarrow u\mbox{ in $L^1_{\rm loc}$ as $n\to\9.$}\end{equation}
 \eqref{e5.7} yields$$\sqrt{u_n b(u_n)}\longrightarrow \sqrt{ub(u)}\mbox{ in }L^2_{\rm loc}$$and so, by Hypothesis (iii),   we have
\begin{equation}\label{e5.8}
E\sqrt{u_nb(u_n)}\longrightarrow E\sqrt{ub(u)}\mbox{ in }L^2_{\rm loc}(\rr^d;\rr^d).\end{equation}
Hence \eqref{e5.6} implies that (selecting a subsequence if necessary) for all balls $B_N$ of radius $N\in\nn$ around zero we have
$$\sup_n\int_{B_N}|\nabla j(u_n)|^2dx<\9$$and$$j(u_n)\to j(u)\mbox{\ \ in }L^2_{\rm loc}\mbox{\ \ as }n\to\9.$$Therefore (again selecting a subsequence, if necessary), for every $N\in\nn$,
$$\nabla j(u_n)\to\nabla j(u)\mbox{ weakly in }L^2(B_N,dx)\mbox{ as }n\to\9.$$
Hence, if we define $\Psi_N$ analogously to $\Psi$, but with the integral over $\rr^d$ replaced by an integral over $B_N$, we conclude that
$$\barr{ll}
\dd\liminf_{n\to\9}\Psi_N(u_n)\ge\!\!\!&\dd\liminf_{n\to\9}\int_{B_N}|\nabla j(u_n)|^2dx-2\int_{B_N}\nabla j(u)\cdot E\sqrt{ub(u)}dx\vsp
&\dd+\int_{B_N}|E|^2ub(u)dx\ge\Psi_N(u).\earr$$
Hence, since $u\in L^1$, we can let $N\to\9$ to get
$$\liminf_{n\to\9}\Psi(u_n)\ge\Psi(u).$$
Now, we consider   the functional
\begin{equation}\label{e5.8az}
\barr{ll}
\Psi_\vp(u)=\!\!\!&\dd\int_{\rr^d}
\left|\frac{\beta'(u)\nabla u}{\sqrt{b^*_\vp(u)+\vp^{2m}}}-E_\vp\sqrt{b^*_\vp(u)+\vp^{2m}}\right|^2dx\vsp
&+\vp^{2m}\dd\int_{\rr^d}E_\vp\cdot\(\frac{\beta'(u)\nabla u}{b^*_\vp(u)+\vp^{2m}}-E_\vp\)dx\vsp
&\dd+\vp\int_{\rr^d}  \beta (u)(\eta'_\vp(u)+\Phi_\vp)dx,\ \ff u\in D(\Psi_\vp)=D(V)\cap H^1,\earr\end{equation}
and
$$\Psi_\vp(u):=\9\mbox{\ \ if }u\in D(V)\setminus H^1.$$
We have

\begin{lemma}\label{l5.2a} For each $\vp>0$, $\Psi_\vp$ is $L^1_{\rm loc}$-lower semicontinuous on every ball in $\calm$. Moreover, for any sequence $\{v_\vp\}\subset D(V)\cap H^1$ such that
$$\sup_{\vp\ge0}\|v_\vp\|<\9,\ \lim_{\vp\to0}v_\vp=v\mbox{ in }L^1_{\rm loc},$$we have
\begin{equation}\label{e5.8a}
\liminf_{\vp\to0}\Psi_\vp(v_\vp)\ge\Psi(v).
\end{equation}
Furthermore, there exists $c\in(0,\9)$ such that, for all $u\in D(V)$, $\vp\in(0,1]$,
\begin{equation}\label{e5.10u}
\Psi_\vp(u)\ge-c(|u|+\|u\|+1).
\end{equation}\end{lemma}

\noindent{\bf Proof.} First of all we note that by the assusmption on $u_\vp$ it follows that $\lim\limits_{\vp\to0} v_\vp=0$ in $L^1$, since $\lim\limits_{|x|\to\9}\Phi(x)=\9.$ We write
 $\Psi_\vp(u)\equiv\Psi^*_\vp(u)+G_\vp(u),$ where
$$\barr{ll}
\Psi^*_\vp(u)=\!\!\!
&\dd\int_{\rr^d}\left|\frac{\beta'(u)\nabla u}{\sqrt{b^*_\vp(u)+\vp^{2m}}}-E_\vp
\sqrt{b^*_\vp(u)+\vp^{2m}}\right|^2dx\vsp
&+\vp^{2m}\dd\int_{\rr^d}E_\vp
\cdot\(\frac{\beta'(u)\nabla u}{b^*_\vp(u)+\vp^{2m}}-E_\vp\)dx,\vsp
G_\vp(u)=\!\!\!
&\vp\dd\int_{\rr^d}\beta(u)(\eta'_\vp(u)+\Phi_\vp)dx.\earr$$
We have, since $\eta'_\vp(\tau)\ge\frac{\gamma_1}{b_0}\,(\log\tau-\vp(1-\tau))$ for $\tau\in(0,1],$
\begin{equation}\label{e5.11u}
\barr{ll}
G_\vp(v_\vp)
 \!\!\!
 &\ge\vp \gamma_1\dd\int_{\{v_\vp\le1\}}
v_\vp\eta'_\vp(v_\vp)dx
\ge\vp\,\dd\frac{\gamma^2_1}{b_0}\int_{\{v_\vp\le1\}}(v_\vp\log v_\vp-\vp v_\vp)dx\vsp
 &\ge\dd-\vp\,\frac{\gamma^2_1}{b_0}
\Bigg[C_\alpha\(\int_{\rr^d}\Phi^{-m}dx\)^{1-\alpha}
\|v_\vp\|^\alpha+\vp\int_{\rr^d}v_\vp\,dx\Bigg],\earr
\end{equation}where we used \eqref{e4.5z}. Hence
 $\liminf\limits_{\vp\to0}G_\vp(v_\vp)\ge0.$
Now, arguing as in the proof of Lemma \ref{l5.1}, we represent $\Psi^*_\vp$ as (see \eqref{e52a})
$$\Psi^*_\vp(u)\!=\!\!\dd\int_{\rr^d}\!
|\nabla j^*_\vp(u){-}E_\vp\sqrt{b^*_\vp(u){+}\vp^{2m}}|^2dx
{+}\vp^{2m}\dd\int_{\rr^d}\!\!E_\vp\cdot\(\dd\frac{\beta'(u)\nabla u}{b^*_\vp(u){+}\vp^{2m}}{-}E_\vp\)dx,$$
where $u\in D(V)\cap H^1$ and 
$$j^*_\vp(r)=\int^r_0\frac{\beta'(s)ds}{\sqrt{b^*_\vp(s)+\vp^{2m}}}.$$
We may assume that $\Psi^*_\vp(v_\vp)\le\nu<\9$, $\ff\vp>0.$ Then, as in \eqref{e52a}, we see that
\begin{equation}\label{e510a}
\barr{ll}
\dd\int_{\rr^d}\frac{|\beta'(v_\vp)|^2|\nabla v_\vp|^2}{b^*_\vp(v_\vp)+\vp^{2m}}\,dx
\!\!\!&\le 2\(\Psi^*_\vp(v_\vp)+\dd\int_{\rr^d}|E_\vp|^2
(b^*_\vp(v_\vp)+2\vp^{2m})dx\)\vsp\qquad
&+2\vp^{2m}\dd\int_{\rr^d}
\frac{|E_\vp|\beta'(v_\vp)|\nabla v_\vp|}{b^*_\vp(v_\vp)+\vp^{2m}}\ dx.\earr
 \end{equation}
Taking into account that
\begin{equation}\label{e5.12u}
\barr{l}
2\vp^{2m}\dd\int_{\rr^d}
\dd\frac{|E_\vp|\beta'(v_\vp)|\nabla v_\vp|}{b^*_\vp(v_\vp)+\vp^{2m}}\,dx\vsp
\qquad\le\dd\frac12\int_{\rr^d}
\dd\frac{|\beta'(v_\vp)|^2|\nabla v_\vp|^2}{b^*_\vp(v_\vp)+\vp^{2m}}\,dx
+{2\vp^{4m}}\int_{\rr^d}
\dd\frac{|E_\vp|^2}{b^*_\vp(v_\vp)+\vp^{2m}}\,dx\vsp
\qquad\le\dd\frac12\int_{\rr^d}
\dd\frac{|\beta'(v_\vp)|^2|\nabla v_\vp|^2}{b^*_\vp(v_\vp)+\vp^{2m}}\,dx
+{2\vp^{2m}} \int_{\rr^d}|E_\vp|^2dx,\earr
\end{equation}
and that $\lim\limits_{\vp\to0}v_\vp=v$ in $L^1$ by our assumption, it follows by \eqref{e3.4x} and \eqref{e510a} that, for some $C>0$ independent of $\vp$,
$$\int_{\rr^d}\frac{|\nabla v_\vp|^2}{b^*_\vp(v_\vp)+\vp^{2m}}\ dx\le C,\ \ff \vp>0,$$and so $\{\nabla j^*_\vp(v_\vp)\}$ is bounded in $L^2$. Then, arguing as in Lemma \ref{l5.1} (see \eqref{e5.7}-\eqref{e5.8}), we get for $\vp\to0$
$$E_\vp \sqrt{b^*_\vp(v_\vp)+\vp^{2m}}\longrightarrow E\sqrt{b(u)u}\mbox{\ \ in }L^2(\rr^d;\rr^d),$$and, therefore,
$$
\dd\liminf_{\vp\to0}\Psi_\vp(v_\vp)\ge
\liminf_{\vp\to0}\Psi^*_\vp(v_\vp)\ge\Psi(v),$$
as claimed. By a similar (even easier) argument, one proves that $\Psi_\vp$ is $L^1_{\rm loc}$-lower semicontinuous on balls in $\calm$. The last part of the assertion is an immediate consequence of \eqref{e5.11u} and \eqref{e5.12u}, which hold for all $u\in D(V)\cap H^1$ replacing $v_\vp$. Hence, the lemma is proved.

We denote by $S_\vp(t)$ the continuous semigroup of contractions on $L^1$ ge\-ne\-rated by the $m$-accretive operator $A_\vp$ defined by \eqref{e3.1}, \eqref{e3.2}, \eqref{e3.7a}, that~is,
\begin{equation}\label{e5.3a}
S_\vp(t)u_0=\lim_{n\to\9}\(I+\frac tn\,A_\vp\)^{-n}u_0,\ \ff t\ge0,\ u_0\in L^1.\end{equation}We note that  by \eqref{e3.3}  it follows, by virtue of the Trotter-Kato theorem for nonlinear semigroups of contractions, that (see \cite{5} and \cite{1}, p.~169)  
\begin{equation}\label{e5.4a}
\lim_{\vp\to0}S_\vp(t)u_0=S(t)u_0,\ \ff u_0\in L^1,\end{equation} 
strongly in $L^1$ uniformly on compact time intervals.

We shall prove first \eqref{e4.6} for $S_\vp(t)$. Namely, one has
\begin{lemma}\label{l5.2} For each $u_0\in L^2\cap D(V)$, we have $S_\vp(\sigma)u_0\in D(\Psi_\vp)$ for $ds$-a.e.  $\sigma\ge0$, and
\begin{equation}\label{e5.9}
V_\vp(S_\vp(t)u_0)+\int^t_s\Psi_\vp(S_\vp(\sigma)u_0)d\sigma\le V_\vp(S_\vp(s)u_0),\   0\le s\le t<\9,\end{equation}
and all three terms are finite.
\end{lemma}

\noindent{\bf Proof.} First, we shall prove   that, for all $\vp>0$,
\begin{equation}\label{e5.10}
V_\vp((I+\lbb A_\vp)^{-1}u_0)+\lbb\Psi_\vp((I+\lbb A_\vp)^{-1}u_0)\le V_\vp(u_0),\  \lbb\in(0,\lbb_0).\end{equation}
We set $u^\lbb_\vp=(I+\lbb A_\vp)^{-1}u_0$ and note that, by \eqref{e3.12c}-\eqref{e3.13a}, we have
\begin{eqnarray}
& u^\lbb_\vp\in H^1(\rr^d),\ \beta(u^\lbb_\vp)\in H^1(\rr^d),\ \ff\lbb\in(0,\lbb_0),\ \vp>0,\label{e5.15prim}\\[2mm]
&V'_\vp(u^\lbb_\vp)=\eta'_\vp(u^\lbb_\vp)+\Phi_\vp\in\partial V_\vp(u^\lbb_\vp),\label{e516a}
\end{eqnarray}
where $\eta'_\vp(u^\lbb_\vp)\in H^1(\rr^d).$
Taking into account that, by Lemma \ref{l3.2},
\begin{equation}\label{e5.15tert}
{\rm div}(\nabla\beta(u^\lbb_\vp)-E_\vp b^*_\vp(u^\lbb_\vp))
=\frac1\lbb\,(u^\lbb_\vp-u_0)+\vp\beta(u^\lbb_\vp)\in\calm\cap L^2,\end{equation}
it follows, since $\Phi_\vp\in L^2$ and ${\rm div}\,E_\vp\in L^2+L^\9$ by \eqref{e340a} and Hypothesis (iii), that
$$\dd\int_{\rr^d}(-\Delta\beta(u^\lbb_\vp)+{\rm div}\,E_\vp b^*_\vp(u^\lbb_\vp))\Phi_\vp\,dx
=-\dd\int_{\rr^d}(\nabla\beta(u^\lbb_\vp)-E_\vp b^*_\vp(u^\lbb_\vp))\cdot E_\vp\,dx.
$$
This yields, by \eqref{e516a},

$$\barr{l}
\<A_\vp(u^\lbb_\vp),V'_\vp(u^\lbb_\vp)\>_2\vsp
\quad=\<-\D(\beta  (u^\lbb_\vp))+\vp\beta  (u^\lbb_\vp)+{\rm div}(E_\vp  b^*_\vp (u^\lbb_\vp)),\eta'_\vp(u^\lbb_\vp)+\Phi_\vp\>_2\vsp
\quad=\dd\int_{\rr^d}(\beta'  (u^\lbb_\vp)\nabla u^\lbb_\vp-E_\vp  b^*_\vp (u^\lbb_\vp))\cdot
\(\dd\frac{\beta' (u^\lbb_\vp)}{b^*_\vp (u^\lbb_\vp)+\vp^{2m}}\nabla u^\lbb_\vp-E_\vp\)dx\vsp
\quad\qquad+\vp\<\beta(u^\lbb_\vp),
\eta'_\vp(u^\lbb_\vp)+\Phi_\vp\>_2\vspace*{3mm}\\
\quad=\dd\int_{\rr^d}\left|\frac{\beta'   (u^\lbb_\vp)\nabla u^\lbb_\vp}{\sqrt{b^*_\vp (u^\lbb_\vp)+\vp^{2m}}}-E_\vp \sqrt{b^*_\vp(u^\lbb_\vp)+\vp^{2m}}\right|^2 dx
+\vp
\<\beta(u^\lbb_\vp),\eta'_\vp(u^\lbb_\vp)
+\Phi_\vp\>_2\vspace*{3mm}\\
\quad\qquad+\vp^{2m}\dd\int_{\rr^d}\(\!E_\vp{\cdot}
 \dd\frac{\beta'(u^\lbb_\vp)\nabla u^\lbb_\vp}{b^*_\vp+\vp^{2m}}-E_\vp\!\)dx
 =\Psi_\vp(u^\lbb_\vp),\, \ff\vp>0,\, \lbb\in(0,\lbb_0).\earr$$This yields \eqref{e5.10}  because, by the convexity of $V_\vp$, we have by \eqref{e516a}
$$V_\vp(u^\lbb_\vp)\le V_\vp(u_0)+\<V'_\vp(u^\lbb_\vp),u^\lbb_\vp-u_0\>_2,\
u^\lbb_\vp-u_0=-\lbb A_\vp(u^\lbb_\vp).$$
  To get \eqref{e5.9}, we shall proceed as in the proof of Theorem 3.4 in \cite{9}. Namely, we set
 $$\barr{r}
 \lbb\delta(\lbb,v)=V_\vp((I+\lbb A_\vp)\1v)+\lbb\Psi_\vp((I+\lbb A_\vp)\1v)-V_\vp(v),\vsp \ff\lbb\in(0,\lbb_0),\ v\in L^2\cap D(V),\earr$$and note that, by \eqref{e5.10}, $\delta(\lbb,u_0)
 \le0,$ $\lbb\in(0,\lbb_0).$ This yields
 $$\barr{c}
 V_\vp((I+\lbb A_\vp)^{-j}u_0)+\lbb\Psi_\vp((I+\lbb A_\vp)^{-j}u_0) -V_\vp((I+\lbb A_\vp)^{-j+1}u_0)\vsp
 =\lbb\delta(\lbb,(I+\lbb A_\vp)^{-j+1}u_0),\ \ff j\in\nn.\earr$$
 Then, summing up from $j=1$ to $j=n$ and taking $\lbb=\frac tn$, we get
\begin{equation}\label{e517a}
\barr{l}
 V_\vp\(\(I+\dd\frac tn\,A_\vp\)^{-n}u_0\)
 +\dd\sum^n_{j=1}\frac tn
 \Psi_\vp\(\(I+\dd\frac tn\,A_\vp\)^{-j}u_0\)\vsp
 \qquad=V_\vp(u_0)+\dd\sum^n_{j=1}\frac tn\ \delta\(\frac tn,\(I+\frac tn\ A_\vp\)^{-(j-1)}u_0\).\earr
 \end{equation}
 Note also that, if $n>\frac t{\lbb_0}$, then
 \begin{equation}\label{e517aa}
\delta\(\frac tn,\(I+\frac tn\ A_\vp\)^{-j}u_0\)\le0,\ \ 1\le j\le n.\end{equation}
   We consider the step function
 $$f_n(\sigma)=\Psi_\vp\(\(I+\frac tn\ A_\vp\)^{-j}u_0\)\mbox{ for }\frac{(j-1)t}n<\sigma\le\frac{jt}n,$$and note that, for each $t>0$,
 $$\sum^n_{j=1}\frac tn\ \Psi_\vp
 \(\(I+\dd\frac tn\,A_\vp\)^{-j}u_0\)=\int^t_0 f_n(\sigma)d\sigma.$$Then, by \eqref{e3.37x}, \eqref{e5.3a} and the $L^1_{\rm loc}$-lower semicontinuity of $\Psi_\vp$ on balls in $\calm$, we conclude, by  the Fatou lemma, which is applicable because of \eqref{e5.10u}, that
\begin{equation}\label{e517aaa}
-\9<\int^t_0\Psi_\vp(S(\sigma)u_0)d\sigma\le
 \liminf_{n\to\9}\int^t_0 f_n(\sigma)d\sigma,
 \end{equation}
 while, by the $L^1_{\rm loc}$-lower semicontinuity of $V_\vp$ on balls in $\calm$, we have
 $$\liminf_{n\to\9}V_\vp
 \(\(I+\dd\frac tn\,A_\vp\)^{-n}u_0\)\ge V_\vp(S_\vp(t)u_0).$$
 Then, by \eqref{e517a}-\eqref{e517aaa}, we get
 $$V_\vp(S_\vp(t)u_0)+\int^t_0\Psi_\vp(S_\vp(\sigma)u_0)d\sigma\le V_\vp(u_0),\ \ff t\ge0.$$In particular, $V_\vp(S_\vp(t)u_0)<\9$ since $V_\vp(u_0)<\9.$ Taking this into account and that $S_\vp(t+s)u_0=S_\vp(t)S_\vp(s)u_0$,
  we get \eqref{e5.9}, as claimed.

\bigskip
\noindent{\bf Proof of Theorem \ref{t4.1} (continued).}  We shall assume \mbox{$u_0\in L^2\cap D_0(V)$.} We want to let $\vp\to0$  in \eqref{e5.9}, where $s=0$.

We note first that  we have
\begin{equation}\label{e5.13a}
\liminf_{\vp\to0}  V_\vp(S_\vp(t)u_0) \ge V(S(t)u_0),\ \ff t\ge0.\end{equation}Here is the argument. We note that, if $v_\vp\!\to\!v$ in $L^1$ as $\vp\!\to\!0$ and \mbox{$\sup\limits_{\vp>0}\|v_\vp\|<\9$,} then $v_\vp(\log v_\vp)^-\to v(\log v)^-$ in $L^1_{\rm loc}$ as $\vp\to0$. Furthermore, for $\delta>0$, and $\alpha \in\left[\frac{m+\delta}{m+\delta+1},1\right)$, by \eqref{e4.5z},
$$\int_{\{\Phi\ge R\}}v_\vp(\log v_\vp)^-dx
\le C_\alpha\ \frac1{R^{\delta(1-\alpha)}}
\(\int\Phi^{-m}dx\)^{1-\alpha}\|v_\vp\|^\alpha,$$
hence
$$\lim_{R\to\9}\sup_{\vp>0}\int_{\{\Phi\ge R\}}v_\vp(\log v_\vp)^-dx=0,$$therefore, $v_\vp(\log v_\vp)^-\to v(\log v)^-$ in $L^1$. Applying this to $v_\vp=S_\vp(t)u_0$, which by \eqref{e5.4a}, \eqref{e3.37x} and \eqref{e5.3a} is justified, and because $\eta_\vp\to\eta$ as $\vp\to0$ locally uniformly on $[0,\9)$ and, because for all $\vp\in(0,1]$, $r\in[0,\9)$,
$$\eta_\vp(r)\ge-\frac{\gamma_1}{b_0}\ (r\wedge1)(\log(r\wedge1)^--2(r\wedge1)),$$we can apply the generalized Fatou lemma to conclude that
$$\liminf_{\vp\to\9}\int_{\rr^d}\eta_\vp(S_\vp(t)u_0)dx\ge\int_{\rr^d}\eta(S(t)u_0)dx,$$and we get \eqref{e5.13a}, as claimed.

By Lemma \ref{l5.2}, \eqref{e3.37x} and \eqref{e5.3a}, we have that $v_\vp=S_\vp(t)u_0$, $\vp>0$, satisfy for $dt$-a.e. $t>0$ the assumptions of Lemma \ref{l5.2a}, hence
$$\liminf_{\vp\to0}\Psi_\vp(S_\vp(t)u_0)\ge
\Psi(S(t)u_0),\ \mbox{a.e. }t>0.$$   Moreover, by Fatou's lemma, which is applicable by \eqref{e5.10u}, it follows that
\begin{equation}\label{e5.25u}
\liminf_{\vp\to0}\int^t_0\Psi_\vp(S_\vp(s)u_0)ds\ge\int^t_0\Psi(S(s)u_0)ds,\ \ff t\ge0.\end{equation}
Because, as mentioned earlier, $V_\vp(u)\to V(u)$ as $\vp\to0$, if $u\in D(V)\cap L^2$, \eqref{e5.13a}, \eqref{e5.25u} and \eqref{e5.9} with $s=0$ imply
\begin{equation}\label{e5.25uu}
V(S(t)u_0)+\int^t_0\Psi(S(\sigma)u_0)d\sigma\le V(u_0),\ \ff u_0\in D(V)\cap L^2,\ t\ge0.\end{equation}
We note that, by \eqref{e2.21a}  and \eqref{e4.3am}, we have
\begin{equation}\label{e5.11}
\barr{lcl}
V(S(t)u_0)&\ge&-C(\|S(t)u_0\|+1)^\alpha\vsp
&\ge&-C(\|u_0\|+
 t|u_0|_1)^\alpha,\ \alpha\in
 \left[{\frac m{m+1}},1\).\earr
\end{equation} 
\n Hence
$$0\le\int^t_0\Psi(S(\sigma)u_0)d\sigma<\9,\ \ff t\ge0,$$which implies that
\begin{equation}\label{e5.25uuu}
S(\sigma)u_0\in D(\Psi)\mbox{\ \ a.e. }\sigma>0.\end{equation}
Now, to extend \eqref{e5.25uu} to all $u_0\in D_0(V),$ take $u^n_0\in D(V)\cap L^2(\subset D_0(V))$ with $u^n_0\le  u_0$ and $u^n_0\to u_0$ as $n\to\9$ in $L^1$. Then, because for all $r\ge0$
$$\eta(r)\ge-\frac{\gamma_0}{b_0}\left[(r\wedge1)(\log (r\wedge1)^-+(r\wedge1))\right],$$arguing as above (using again \eqref{e4.5z}), we conclude the monotone convergence applies to get
$$\lim_{n\to\9} V(u^n_0)=V(u_0)$$and the generalized Fatou lemma applies to get eventually \eqref{e5.25uu} and \eqref{e5.25uuu} for all $u_0\in D_0(V).$ Since $S(t)u_0\in D_0(V)$, if $u_0\in D_0(V)$, the first part including \eqref{e4.6} follows.

To prove \eqref{e4.11i}, we note that since $\alpha<1,$  by \eqref{e4.6} and \eqref{e5.11},  we have
\begin{equation}\label{e5.30'}
\barr{lcl}
0&=&\dd\lim_{t\to\9}\frac1t\int^t_0\Psi(S(\sigma)u_0)d\sigma \ge \dd\lim_{t\to\9}\frac1t\int^t_n\inf_{r\ge n}\Psi(S(r)u_0)d\sigma\vsp
&=&\dd\inf_{r\ge n}\Psi(S(r)u_0)\mbox{\ \ for all }n\in\nn.\earr\end{equation}
Hence, there exists $t_n\to\9$ such that
\begin{equation}\label{e5.31i}
\lim_{n\to\9}\Psi(S(t_n)u_0)=0.
\end{equation}Furthermore, we obtain by Lemma \ref{l5.1} the first inequality in \eqref{e5.30'}, \eqref{e2.15} and \eqref{e2.24'}  that
$$\sup_{t\ge0}|S(t)u_0|_1+\limsup_{t\to\9}
\frac1t\int^t_0|\nabla(\sqrt{S(s)u_0})|_2ds<\9.$$
Hence, similarly as above (selecting a subsequence of $(t_n)$, if necessary),
\begin{equation}\label{e5.26u}
\sup_n\|\sqrt{S(t_n)u_0}\|_{W^{1,2}(\rr^d)}<\9.\end{equation}
So, by the Rellich-Kondrachov theorem (see, e.g., \cite{5}, p.~284), the set
$$\{S(t_n)u_0\mid n\in\nn\}$$ is relatively compact in $L^1_{\rm loc}$. Hence, along a subsequence
$\{t_{n'}\}\to\9$, we~have 
\begin{equation}\label{e5.18}
\lim S(t_{n'})u_0=u_\9\mbox{ in }L^1_{\rm loc}
\end{equation}
for some $u_\9\in L^1.$ Since $\Psi$ is $L^1_{\rm loc}$-lower semicontinuous on $L^1$-balls by Lemma \ref{l5.1}, this together with \eqref{e5.31i}  implies  that $u_\9\in D(\Psi)$ and $\Psi(u_\9)=0$.

If $u_\9\in D(\Psi)$, such that $\Psi(u_\9)=0$, then
\begin{equation}\label{e6.12u}
\frac{\beta'(u_\9)\nabla u_\9}{\sqrt{u_\9 b(u_\9)}}=E\sqrt{u_\9 b(u_\9)},\mbox{ a.e. in }\rr^d.\end{equation}
Let us prove now that either $u_\9\equiv0$ or $u=u_\9>0$, a.e. in $\rr^d$. To this end, we consider the solution $y=y(t,x)$ to the system
$$\barr{l}
y'_i(t)=\wt D_i (y_i(t)),\ t\ge0,\ i=1,...,d,\vsp
y_i(0)=x_i,\earr$$
where $\wt D_i\in  C^1(\rr),$ $i=1,....,d,$  is an arbitrary vector field on $\rr$ of at most linear growth, and   $y(t)=\{y_i(t)\}^d_{i=1},$ $x=\{x_i\}^d_{i=1}.$ If $j$ is defined by \eqref{e5.4c}, we have
$$\barr{ll}
\dd\frac d{dt}\,j(u(y(t,x)))\!\!\!
&=j_u(u(y(t,x)))\nabla u(y(t,x))\cdot\dd\frac d{dt}\,
y(t,x)\vsp
&=\dd\frac{\beta'(u(y(t,x)))}
{\sqrt{b(u(y(t,x)))u(y(t,x))}}\,\nabla u(y(t,x)){\cdot}\cald (y(t,x)), \ff t\ge0,\earr$$
where
$\cald(y)=(\wt D_{i}(y_i))^d_{i=1}$. Let $E=\{E_i\}^d_{i=1}$. Then, by \eqref{e6.12u},
$$\frac d{dt}\,j\(u(y(t,x))\)=
\sum^d_{i=1}\wt D_i (y_i(t,x))E_i(u(y(t,x)))
\(u(y(t,x))  b(u(y(t,x)))\)^{\frac12}.$$
We note that
$$C_2j(r)\le\sqrt{rb(r)}\le C_1j(r), \ff r\ge0,$$where $C_1,C_2>0$. We set $\alpha(t,x)=(u(y(t,x))b(u(y(t,x))))^{\frac12}(j(u(y(t,x))))\1$. Then $\alpha\in L^\9((0,\9)\times\rr^d)$ and
$$\frac d{dt}\,j(u(y(t,x)))=\alpha(t,x)\sum^d_{i=1}
\wt D_i(y_i(t,x))E_i(u(y(t,x)))j(u(y(t,x))),\ \ff t\ge0.$$Hence
$$j(u(y(t,x)))=j(u(x))\exp\(\int^t_0\alpha(s,x)\cald(e^{\cald s}x)\cdot E(u(e^{\cald s}x))\), \ff t\ge0,\ x\in\rr^d,$$and, therefore, 
$$j(u(x))=j(u(e^{\cald t}x))\exp\(-\int^t_0\alpha(s,x)\cald(e^{\cald s}x)\cdot E(u(e^{\cald s}x))\),$$
where $e^{\cald t}$ is the flow generated by $\cald$.
Since $\cald$ is an arbitrary vector field on $\rr^d$, it follows that, for fixed $x$ and $t$, $\{e^{\cald t}x\}$ covers all $\rr^d.$ We infer that, if $u\not\equiv0$, then $j(u(x))>0$, $\ff x\in\rr^d$, and this implies that $u=u_\9>0$, a.e. on $\rr^d$.
For such a $u_\9$, this yields, because $\Psi(u_\9)=0,$
\begin{equation}\label{e6.13u}
\nabla(g(u_\9)+\Phi)=0,\mbox{ a.e. in }\rr^d,\end{equation}where
$$g(r)=\int^r_1\frac{\beta'(s)}{sb(s)}\ ds,\ \ff r>0.$$By \eqref{e6.13u}, we see that $g(u_\9)+\Phi=\mu$ for some $\mu\in\rr,$ in $\rr^d$ and, since $g$ is strictly monotone, we have
\begin{equation}\label{e6.14u}
u_\9(x)=g^{-1}(-\Phi(x)+\mu),\ \  x\in\rr^d.
\end{equation}

\section{The asymptotic behaviour in $L^1$}\label{s6}
\setcounter{equation}{0}

\begin{theorem}\label{t6.1} Assume that Hypotheses {\rm(i)-(vi)}    hold and let \mbox{$u_0\!\in\! D_0(V){\setminus}\{0\}$.} Set
$$\wt\oo(u_0)=\left\{\lim\limits_{n\to\9}S(t_n)u_0\mbox{\ \ in }L^1,\ \{t_n\}\to\9\right\}.$$  Then
\begin{equation}\label{e6.1}
 \oo(u_0)=\wt\oo(u_0)=\{u_\9\},\end{equation}
 and $u_\9 > 0$, a.e. on $\rr^d$. Furthermore, $u_\9 \in  D_0(V)\cap D(\Psi),$ $ \Psi(u_\9)=0,$ $  S(t)u_\9=u_\9$ for $t\ge0$, $|u_\9|_1=|u_0|_1$, and  it is given by
 \begin{equation}
 u_\9(x)=g\1(-\Phi(x)+\mu),\ \ff x\in\rr^d, \label{e6.2}\end{equation}
where $\mu$ is the unique number in $\rr$ such that
 \begin{equation}
 \dd\int_{\rr^d}g\1(-\Phi(x)+\mu)dx=\int_{\rr^d}u_0\,dx,
 \label{e6.3}
\end{equation}where
$$g(r)=\int^r_1\frac{\beta'(s)}{sb(s)}\,ds,\ r>0.$$
In particular, for all $u_0\in D_0(V)$ with the same $L^1$-norm, the sets in \eqref{e6.1} coincide, and thus $u_\9$ is the only element in $D_0(V)$ with given $L^1$-norm such that $S(t)u_\9=u_\9$ for all $t\ge0$.
\end{theorem}
\noindent{\bf  Proof.}  Let us first   prove the following version of Proposition \ref{p2.3}.

\begin{lemma}\label{l6.2} Under Hypotheses {\rm(i)-(vi)}, we have, for all $u_0\in\calm_+,$
 \begin{eqnarray}
 \|(I+\lbb A)^{-1}u_0\|&\le&\|u_0\|,\ \ff \lbb\in(0,\lbb_0),\label{e6.3aa}\\
\|S(t)u_0\|&\le&\|u_0\|,\ \ff  t\ge0.\label{e6.4}
 \end{eqnarray}
\end{lemma}

\noindent{\bf Proof.} We may assume that by approximation  $u_0\in\calm_+\cap L^2$. Arguing as in the proof of Lemma \ref{l3.2} and taking into account that $u_\vp\ge0$,  we get by \eqref{e3.36av}-\eqref{e3.37a},
\begin{equation}\label{e6.8}
\hspace*{-4mm}\barr{r}
\dd\int_{\rr^d}\!\!\!u_\vp\vf_\nu dx
 \!\le\!-\lbb\!\!\dd\int_{\rr^d}\!\!\!
((b^*_\vp(u_\vp)|\nabla\Phi_\vp|^2{+}
\nabla\Phi_\vp{\cdot}
\nabla\beta(u_\vp))(1{-}\nu\Phi_\vp)
\exp(-\nu\Phi_\vp))dx\\
 +\dd\int_{\rr^d}u_0\vf_\nu\,dx.\earr\hspace*{-4mm}\end{equation}
 Since, by \eqref{e33a} and Hypotheses (iii), (iv), we have that $|\nabla\Phi_\vp|\in L^2$ and $\beta(u_\vp)\in H^1$, we may pass to the limit $\nu\to0$ in \eqref{e6.8} to find after integrating by parts using Hypothesis (v) that
 \begin{equation}
 \label{e6.9n}
 \barr{ll}
 \dd \int_\rrd u_\vp\Phi_\vp\,dx\le\!\!\!&\dd \lbb
 \int_\rrd\(-b_0\,.\frac{u_\vp}{1+\vp|u_\vp|}\,|\nabla\Phi_\vp|^2+\Delta\Phi_\vp\beta(u_\vp)\)dx\vsp\qquad\qquad
 &+\dd\int_\rrd u_0\Phi_\vp\,dx.\earr
 \end{equation}
 We note that  integrating by parts is justified here, since $\beta(u_\vp)\in L^1\cap L^2$ and $\Delta\Phi_\vp\in L^2+L^\9$ because of \eqref{e340a} and Hypothesis (iii). Now, we want to let $\vp\to0$ (along a subsequence) in \eqref{e6.9n}. To this end, we note that, since by Hypothesis (iii)
 $\Delta\Phi=f_2+f_\9$ for some $f_2\in L^2$, $f_\9\in L^\9$, it follows by \eqref{e340a}, \eqref{e3.45prim} that
 $$\Delta\Phi_\vp=g_\vp(f_2+f_\9)+\vp h_\vp|D|^2,$$where $g_\vp,h_\vp:\rrd\to\rr$ such that $g_\vp\to1$, a.e.  as $\vp\to0$, with $|g_\vp|\le m+1$ and $|h_\vp|\le m(m+3).$ Since $\beta(u_\vp)\to\beta(u)$ in $L^1$ by Lemma \ref{l3.2} (a) and also weakly in $L^2$ by \eqref{e3.18} and since $|\nabla\Phi_\vp|^2\to|\nabla\Phi|^2$, a.e. as $\vp\to0$, by \eqref{e33a}, by virtue of Fatou's lemma we can pass to the limit $\vp\to0$ (along a subsequence) in \eqref{e6.9n} to obtain
 $$\|u\|\le\lbb\int_\rrd(-b_0|\nabla\Phi|^2u+\Delta\Phi\beta(u))dx+\|u_0\|,$$where $u=J_\lbb u_0=(I+\lbb A)^{-1}u_0$ is as in \eqref{e3.36x}. By Hypothesis (vi), this implies \eqref{e6.3aa}, which in turn implies \eqref{e6.4} by the same argument as in the proof of Proposition \ref{p2.3}.

As a consequence of Lemma \ref{l6.2}, inequality \eqref{e2.21a} holds, hence we can apply Theorem \ref{t4.1} below. Hence, by \eqref{e4.3am} and \eqref{e6.4}, we have, for all $t\ge0$,
$$V(S(t)u_0)\ge- C(\|S(t)u_0\|+1)^\alpha\ge-C(\|u_0\|+1)^\alpha,$$hence, by \eqref{e4.6},
\begin{equation}\label{e6.7u}
\int^\9_0\Psi(S(\sigma)u_0)d\sigma<\9.\end{equation}
This implies that
\begin{equation}\label{e6.8u}
\oo(u_0)\subset\{u\in D(\Psi);\ \Psi(u)=0\}.\end{equation}
To prove this, we shall use a modification of the argument from the proof of Theorem 4.1 in \cite{9}.

Let $u_\9\in \oo(u_0)$ and $\{t_n\}\to\9$ such that
$S(t_n)u_0\to u_\9\mbox{ in }L^1_{\rm loc}.$
Assume that $\Psi(u_\9)>\delta>0$ and argue from  this to a contradiction. This  implies that there is a bounded open subset $\calo$ of $\rr^d$ such that
\begin{equation}\label{e6.9u}
\Psi_\calo(u_\9)>\frac\delta2>0,\end{equation}
where $\Psi_{\calo}$ is the integral for \eqref{e4.4} restricted to $\calo$.
Since $\Psi_\calo$ is lower semi\-con\-ti\-nuous in $L^1$, it follows by \eqref{e6.9u} that there is a $\mu=\mu(\delta)>0$ such that
\begin{equation}\label{e6.10u}
\Psi_\calo(u)\ge\frac\delta4\mbox{\ \ if }|u_\9-u|_{1}\le\mu.\end{equation}
Since $S(t),\ t>0$, is a semigroup of contractions, we have
\begin{equation}\label{e6.11u}
|S(t)u_0-S(s)u_0|_{1}\le\nu(|t-s|),\ \ff s,t\ge0,\end{equation}
where $\nu(r):=\sup\{|S(s)u_0-u_0|_{1}:0\le s\le r\},\ r>0.$ Clearly, $\nu(r)\to0$ as $r\to0.$
By \eqref{e6.11u}, we have
$$|S(t)u_0-u_\9|_{1}\le |S(t)u_0-S(t_n)u_0|_{1}
+|S(t_n)u_0-u_\9|_{1}\le\mu,$$
for $|t-t_n|\le\nu^{-1}\(\frac\mu2\),\ n\ge N(\mu)$, where $\nu\1$ is the inverse function of $\nu$.
By~\eqref{e6.10u}, this yields
$$\Psi_\calo(S(t)u_0)\ge\frac\delta4\mbox{ for }|t-t_n|\le\nu^{-1}\(\frac\mu2\),$$and $n\ge N(\mu).$
But this contradicts \eqref{e6.7u}. 

\eqref{e6.8u} and Theorem \ref{t4.1} imply \eqref{e6.2}. By \eqref{e6.4}, we also have
$$\lim_{R\to\9}\ \sup_{t\ge0}\ \int_{\{\Phi\ge R\}} S(t)u_0\ dx=0,$$which implies that the orbit $\{S(t)u_0,\ t\ge0\}$ is compact in $L^1$,  $\oo(u_0)=\wt\oo(u_0)$ and that $|u_\9|_1=|u_0|_1$ by \eqref{e2.11} and \eqref{e2.12}.

Hence \eqref{e6.3} follows and thus \eqref{e6.1} also holds. By Fatou's lemma, it follows that $u_\9\in D(V)$ and, by \eqref{e6.14u}, \eqref{e4.8prim} and the $L^1_{\rm loc}$-lower semicontinuity of $V$ on balls in $\calm$, we conclude that $u_\9\in D_0(V)$. Now, let us check that $S(t)u_\9=u_\9,$ for $t\ge0$.
So, let $t_n\to\9$, such that $\lim_{n\to\9}S(t_n)u_0=u_\9.$
Then, for all $t>0$, by the semigroup property and the $L^1$-continuity of $S(t)$,
$$S(t)u_\9=
\lim_{n\to\9}S(t+t_n)u_0\in\wt\oo(u_0)=\{u_\9\}.$$The last part of the assertion is obvious by \eqref{e6.3}.

\begin{corollary}\label{c6.2} Let $u_\9$ be as in Theorem {\rm\ref{t6.1}.}
	Then
	$$|u_\9|_\9\le\max
	\(1,e^{\frac{|b|_\9}{\gamma}\,(\mu-1)}\),$$where $\mu\in\rr$ is as in \eqref{e6.2}.  \end{corollary}
 \n{\bf Proof.} For $g$ as above, we have that $g$ is strictly increasing and \mbox{$g:(0,\9)\to\rr$} is bijective. Furthermore, by \eqref{e4.1av}, we have, for $r\in(0,\9)$,
 $$\frac{\gamma_1}{b_0}\ {\bf1}_{(0,1]}(r)\log r+
\frac\gamma{|b|_\9}\ {\bf1}_{(1,\9)}(r)\log r\le g(r).$$Hence, replacing $r$ by $e^{\frac{b_0}{\gamma_1}\,r},$ $r\le0$, we get
$$g\1(r)\le e^{\frac{b_0}{\gamma_1}\,r},\ r\in(-\9,0],$$and, replacing $r$ by $e^{\frac{|b|_\9}{\gamma}\,r}$, $r\in(0,\9)$, we obtain
$$g\1(r)\le e^{\frac{|b|_\9}{\gamma}\,r},\ r\in(0,\9).$$
This implies, by \eqref{e6.2}, for all $x\in\rr^d,$
$$\barr{ll}
(0<)u_\9(x)\!\!\!&=g\1(\mu{-}\Phi(x))
\le{\bf1}_{\{\mu\le\Phi\}}(x)e^{\frac{b_0}{\gamma_1}\,(\mu-\Phi(x))}{+}{\bf1}_{\{\mu>\Phi\}}(x)e^{\frac{|b|_\9}\gamma\,(\mu-\Phi(x))}\vsp
&\le\max\(1,e^{\frac{|b|_\9}\gamma\,(\mu-1)}\),\earr$$since $\Phi\ge1.$

  We show now that Theorem \ref{t6.1} implies  the uniqueness of solutions $u^*\in\calm\cap\calp\cap\{V<\9\}$ of the stationary version of \eqref{e1.1}, that is, to the equation
  \begin{equation}\label{e6.12}
  -\Delta\beta(u^*)+{\rm div}(Db(u^*)u^*)=0\mbox{ in }\cald'(\rr^d).\end{equation}
  We note that the set of all $u^*\in L^1(\rr^d)$ satisfying \eqref{e6.12} is just $A^{-1}_0(\{0\}).$

  \begin{theorem}\label{t6.2} Under Hypotheses {\rm(i)-(vi)}, there is a unique solution   $u^*$ to equation \eqref{e6.12}, such that $u^*\in L^1\cap L^\9$. In addition, $u^*\in\calm\cap\calp\cap\{V<\9\}$.
  \end{theorem}

  \n{\bf Proof.} By Theorem \ref{t6.1} and Corollary \ref{c6.2}, it follows that $u_\9$ is a solution to \eqref{e6.12}, which is in $\calm\cap\calp\cap\{V<\9\}\cap L^\9$. So it only remains to  prove  the uniqueness. But this follows from Theorem 2.1 in \cite{7a}.

\begin{theorem}
	\label{t6.4}
Let $X^i(t),\ t\ge0,$ $i=1,2,$ be two stationary nonlinear distorted Brownian motions, i.e., both satisfy \eqref{e1.4a} with  $(\calf_t^i)$-Wiener processes $W^{i}(t)$, $t\ge0$, on probability spaces $(\Omega^i,\calf^i,\mathbb{P}^i)$ equipped with normal filtrations $\calf^i_t$, $t\ge0$, with
$$\mathbb{P}^i\circ(X^i(t))\1=u^i_\9\,dx,$$ and $u(t,x)$ in \eqref{e1.4a} replaced by $u^i_\9(x)$ for $i=1,2,$ respectively. Assume that $u^i_\9\in\calm\cap\{V<\9\}\cap L^\9$, $i=1,2.$ Then
	 $$\mathbb{P}^1\circ(X^1)\1=\mathbb{P}^2\circ(X^2)\1,$$i.e., we have uniqueness in law of stationary nonlinear distorted  Brownian motions with stationary measures in $\calm\cap\{V<\9\}\cap L^\9.$
\end{theorem}

\n{\bf Proof.} By It\^o's formula, both $u^1_\9$ and $u^2_\9$  satisfy \eqref{e6.12}. Hence, by Theorem \ref{t6.2}, we have $u^1_\9=u^2_\9=u_\9$. Fix $T>0$ and let
$$\Phi(r):=\frac{\beta(r)}r,\ r\in\rr.$$Then Theorem 3.1 in \cite{6a}  implies that, for each $s\in[0,T]$ and each $v_0\in L^1\cap L^\9$, there is at most one solution $v=v(t,x)$, $t\in[s,T]$, to
$$\barr{l}
v_t-\Delta(\Phi(u_\9)v)+{\rm div}(Eb(u_\9)v)=0\mbox{ in }\cald'((0,T)\times\rr^d,\vsp
v(0,\cdot)=v_0,\earr$$such that $v\in L^\9((s,
T)\times\rr^d)$ and $t\mapsto v(t,x)dx,$ $t\in[s,T]$ is narrowly continuous. But $u_\9,$ the time marginal law of $X^i$ under $\mathbb{P}^i$, $i=1,2,$ is such a solution with $v_0=u_\9$, since $u_\9\in L^\9$ by Corollary \ref{c6.2}. Hence, Lemma 2.12 in \cite{14prim} implies the assertion, since by It\^o's formula $\mathbb{P}^i\circ(X^i)\1,$ $i=1,2,$ both satisfy the martingale problem for the Kolmogorov operator
$$L_{u_\9}=\Phi(u_\9)\Delta+b(u_\9)E\cdot\nabla.$$
\begin{remark}
	\label{r6.5}
	\rm By \cite{4}, a stationary nonlinear distorted Brownian motion as above always exists under the assumptions in this section. Furthermore, we recall that for $u\in\calm_+$ by definition of $V$ we have $u\in\{V<\9\}$ if and only if $u\log u\in L^1.$
\end{remark}

\section*{Appendix}

Let $\alpha=\frac{b_0}{\gamma_1},\ \delta=\exp\(-\frac{d+2}{2d}\)$ and $\eta,\mu\in(0,\9)$ to be chosen (large enough) later.
Let $h:[\delta,\9)\to\rr$ be the solution to the following ODE:

$$h'(r)+\frac{d-1}r\,h(r)-\alpha h^2(r)=0,\ r\in(\delta,\9),\eqno\rm(A.1)$$
$$h(\delta)=\delta(2\log\delta+1)-\eta\ 
\eqno\rm(A.2)$$
As we shall see below, it is easy to solve (A.1) explicitly. The solution has the following properties:
 (h.1)  $h$ is bounded;
	 (h.2)  $h$ is negative, $|h(r)|\le C|r|(1+\log|r|)\1$, and there exist $C\in(0,\9)$ and $\wt\eta\in(0,\eta)$ such that
	$\int^r_\delta h(s)ds\ge-C-\wt\eta(r-\delta),\ \ r\in[\delta,\9).$

Now, define as in \eqref{e1.3prim} and \eqref{e1.2a}
$$\vf(r)=\delta^2\log\delta-\eta\delta+\int^r_\delta h(s)ds,\ r\in[\delta,\9),\eqno\rm(A.3)$$ 
$$\Phi(x)=\left\{\barr{ll}
|x|^2\log|x|+\mu,&\mbox{ for }|x|\le\delta,\vsp
\vf(|x|)+\eta|x|+\mu&\mbox{ for }|x|>\delta.\earr\right.\eqno\rm(A.4)$$
Then $\Phi\in C(\rr^d)\cap W^{1,1}_{\rm loc}(\rr^d)$ and by (h.2) for large enough $\mu>0$ and some $\vp>0$,
$\Phi(x)\ge1+\vp|x|\mbox{ for }|x|>\delta.$ Furthermore,
$$\nabla\Phi(x)=\left\{\barr{ll}
x(2\log|x|+1)&\mbox{ for }|x|\le\delta,\\
(h(|x|)+\eta)\,\dd\frac x{|x|}&\mbox{ for }|x|>\delta.\earr\right.
\eqno\rm(A.5)$$
By (A.2)  and (h.1), it follows that $E=-\nabla\Phi\in C_b(\rrd;\rrd).$ Since $h'$ is  bounded, it follows that $\nabla\Phi\in W^{1,1}_{\rm loc}(\rrd;\rrd)$
$$
\Delta\Phi(x)=\left\{\barr{ll}
2d\log|x|+d+2&\mbox{ for }|x|\le\delta,\\
h'(|x|)+\dd\frac {d-1}{|x|}(h(|x|)+\eta)&\mbox{ for }|x|>\delta.\earr\right.\eqno\rm(A.6)$$
 Hence, $\Phi$ satisfies both conditions (iii)  and (iv). It remains to show \eqref{e1.2a}. To this end, we first note that
 $\Delta\Phi(x)\le0\le\alpha|\nabla\Phi(x)|^2\mbox{ for }|x|\le\delta.$ Furthermore, for $|x|\ge\delta$, by (A.1) and (A.6),
$$\barr{lcl}
\Delta\Phi(x)&=&\alpha h^2(|x|)+\dd\frac{d-1}{|x|}\,\eta\\
&=&\alpha|\nabla\Phi(x)|^2+\eta\(\dd\frac{d-1}{|x|}-\alpha(2h(|x|)+\eta)\)
\le\alpha|\nabla\Phi(x)|^2,\earr$$by (h.1) and (h.2), if we choose $\eta>0$ large enough. Hence, $\Phi$ satisfies condition (vi). It remains to solve (A.1), (A.2) and prove that (h.1) and (h.2) hold. This is elementary, but we include it for the convenience of the reader.

  Let $I:=[\delta,\inf\{r>\delta\mid h(r)=0\}] $
and $h:I\to\rr$ be such that (A.1), (A.2) hold.
Setting $g:=\frac1h,$ we see that
$$g'(r)-\frac{d-1}r\,g(r)=-\alpha,\ r\in I,\ \
g(\delta)=-\(\frac{2\delta}d+\eta\)^{-1}.\eqno\rm(A.7)$$
We can rewrite (A.7) equivalently as
 $(r^{1-d}g(r))'=-\alpha r^{1-d},\ r\in I.$ Hence,
$$g(r)=\left\{\barr{ll}
r^{d-1}\left[\delta^{1-d}g(\delta)-\dd\frac\alpha{2-d}\,(r^{2-d}-\delta^{2-d})\right],&\mbox{ if }d\ne2,\vsp
r[\delta^{-1}g(\delta)-\alpha(\log r-\log\delta)],&\mbox{ if }d=2,\earr\right.$$which implies that $I=[\delta,\9)$ and that, for $r\ge\delta,$
$$
\hspace*{-4mm}h(r)=\left\{\barr{ll}
-r^{-1}\left[\(\delta^{-1}\(\dd\frac{2\delta}d+\eta\)^{-1}+\frac\alpha{d-2}\)
\(\frac r\delta\)^{d-2}-\frac\alpha{d-2}\right]^{-1},
\hspace*{-3mm}&\mbox{ if }d\ne2,\vsp
-r^{-1}\left[\delta^{-1}\(\delta+\eta\)^{-1}+\alpha\log\frac r\delta\right]^{-1},&\mbox{ if }d=2.
\earr\right.\hspace*{-4mm}\eqno\rm(A.8)$$
So, $h$ is negative and (h.1) holds, since $|h(r)|\le\frac{2\delta}d+\eta,$ $r\in[\delta,\9)$. Now, we show (h.2) for $d=1$, $d=2$, $d\ge3$, separately.

\bigskip\n{\bf Case $d=1$.} In this case with $g(\delta)$ as defined in (A.7), we have, for $r\in[\delta,\9)$,
 $h(r)=-[|g(\delta)|+\alpha(r-\delta)]^{-1} ,$  and hence, for $K\in(1,\9)$, 
$$\dd\int^r_\delta h(s)ds=
-\dd\frac1\alpha\log\(1+\frac\alpha{|g(\delta)|}\,(r-\delta)\)
\ge-\dd\frac1\alpha\log K{-}K\1|g(\delta)|\1(r{-}\delta)$$and so (h.2) follows for $K$ large enough.

\medskip\n{\bf Case $d=2$.} In this case  we have, for $r\in[\delta,\9)$ and $K\in(1,\9)$,
$$\dd\int^r_\delta h(s)ds=
-\dd\frac1\alpha\log\(1+\frac{\delta\alpha}{|g(\delta)|}\,\,\log\frac r\delta\)
\ge-|g(\delta)|\1K\1(r-\delta)$$and (h.2) follows for $K$ large enough.\medskip

 \n{\bf Case $d=3$.} In this case  we have, for $r\in[\delta,\9)$,
 $|h(r)|\le\(\frac r\delta\)^{1-d}|g(\delta)|\1,$ hence, for $K\in(1,\9)$,
$$\barr{lcl}
\dd\int^r_\delta|h(s)|ds&\le&
(K-1)\delta|g(\delta)|\1+|g(\delta)|\1\delta^{d-1}
\dd\int^{\max(r,K\delta)}_{k\delta}s^{1-d}ds\\
&\le&g(\delta)\1((K-1)\delta+K^{1-\delta}(r-\delta)),\earr$$and (h.2) follows for $K$ large enough.

\bk\n{\bf Acknowledgements.} This work    was supported by the DFG through CRC 1283 and by UEFISCDI (Romania) through PN-III-ID-PCE-2021-3. The authors are indebted to the anonymous referee for carefully rea\-ding this work and for very useful  suggestions.


\begin{thebibliography}{nn}
	
	
 	
 	
\bibitem{15a} Arnold, A., Markowich, P., Toscani, G., Unterreiter, A., On convex Sobolev inequalities and the rate of convergence to equilibrium for Fokker-Planck type equations, {\it Comm. Partial Differential Equations}, vol. 26 (2001), 43-100.\vspace*{-3mm}

\bibitem{1a} Bakry, D., Gentil, I., Ledoux, M., {\it Analysis and geometry of Markov dif\-fu\-sion operators}, Springer, 2014, xx+552 pp. ISBN: \mbox{378-3-319-00226-2.}\vspace*{-3mm}

\bibitem{1} Barbu, V., {\it Nonlinear Differential Equations of Monotone Type in Banach Spaces}, Springer 2010, New York, Dordrecht, Heidelberg, London.\vspace*{-3mm}

\bibitem{2}  Barbu, V., Generalized solutions to nonlinear Fokker-Planck equations, {\it J. Diff. Equations}, 261 (2016), 2446-2471.\vspace*{-3mm}
	
\bibitem{3}  Barbu, V.,  R\"ockner, M., Probabilistic representation for solutions to nonlinear Fokker-Planck equations, {\it SIAM J. Math. Anal.}, 50 (2018), 2588-2607.\vspace*{-3mm} 


\bibitem{4} Barbu, V.,  R\"ockner, M., From nonlinear Fokker-Planck equations to solutions of distribution dependent SDE, {\it Annals of Probability}, 48 (4) (2020), 1902-1920.\vspace*{-3mm}


\bibitem{6a} Barbu, V.,  R\"ockner, M., Solutions for nonlinear Fokker-Planck equations with measures as initial data and  McKean-Vlasov equations, {\it J. Funct. Anal.}, 280 (7) (2021), 1-35.\vspace*{-3mm}

\bibitem{7a} Barbu, V., R\"ockner, M., Uniqueness for nonlinear Fokker-Planck equations and weak uniqueness for McKean-Vlasov SDEs, {\it Stoch. PDEs; Anal. Computation}, 9 (4) (2021).\vspace*{-3mm}

\bibitem{8ab} Barbu, V., R\"ockner, M., Corrections to: Uniqueness for nonlinear Fokker--Planck equations and weak uniqueness for McKean--Vlasov SDEs, {\it Stoch. PDEs; Anal. Computation,} (2022).

\bibitem{8aa} Barbu, V., R\"ockner, M., The existence and uniqueness of nonlinear Fokker--Planck flows (to appear).

\bibitem{0} Bogachev, V.I., Krylov, N.V., R\"ockner, M., Shaposhnikov, S.V., {\it Fokker-Planck-Kolmogorov equations,} Mathematical Surveys and Monographs, 207,   American Mathematical Sociedy, Providence, R.I., 2015, xii+479 pp. ISBN: 978-1-4704-2558-6.\vspace*{-3mm}

\bibitem{0prim} Bogachev, V.I., R\"ockner, M., Shaposhnikov, S.V., Convergence in va\-ria\-tion of solutions of nonlinear Fokker-Planck-Kolmogorov equations to stationary measures, {\it J.~Funct. Anal.}, 276  (12) (2019), 3681-3713.\vspace*{-3mm}


\bibitem{4a} Brezis, H., {\it Functional Analysis Sobolev Spaces and Partial Differential Equations}, Springer 2010, New York, Dordrecht, Heidelberg, London.\vspace*{-3mm}



\bibitem{5} Brezis, H., Pazy, A., Convergence and approximation of semigroups of nonlinear operators in Banach spaces, {\it J. Funct. Anal.}, 9 (1972), 63-74.\vspace*{-3mm}

\bibitem{13a} Carillo, J.A., J\"ungel, A., Markowich, P.A., Toscani, G., Unterreiter, A., Entropy dissipation methods for degenerate parabolic problems and generalized Sobolev inequalities, {\it Monatsh. Math.}, 133 (2001), 1-82.\vspace*{-3mm}

\bibitem{9prim} Carillo, J.A., Toscani, G., Asymptotic $L^1$-decay of solutions of the porous media equation to self-similarity, {\it Indiana Univ. Math. J.,} 49 (1) (2000), 113-142.\vspace*{-3mm}


\bibitem{6} Chavanis, P.H., Generalized stochastic Fokker-Planck equations, {\it Entropy}, 2015, 3205-3252.\vspace*{-3mm}


\bibitem{8a} Chen, G.Q.,   Perthame, B., Well posedness for nonisotropic degenerate parabolic hyperbolic equations, {\it Ann. Institute H. Poincar\'e}, 4 (2003), 645-668.\vspace*{-3mm}

\bibitem{7c}  Crandall, M.G., The semigroup approach to first order quasilinear equations in several space variables, {\it Israel J. Math.,} 10 (1972), 108-132.\vspace*{-3mm}


\bibitem{8b} Dafermos, C., Slemrod, M., Asymptotic behavior of nonlinear contractions semigroups, {\it J. Funct. Anal.}, 13 (1973), 97-100.\vspace*{-3mm}


\bibitem{14primm} Eberle, E., Guillin, A., Zimmer, R., Quantitative Harris-type theorems for diffusions and McKean-Vlasov processes, {\it Trans. Amer. Math. Soc.}, 371 (10) (2019), 7135-7173.\vspace*{-3mm}

\bibitem{7} Frank, T.D., Generalized Fokker-Planck equations derived from ge\-ne\-ra\-lized linear nonequilibrium thermodynamics, {\it Physica} A, 310  (2002), 397-412.\vspace*{-3mm}

\bibitem{8} Frank, T.D., {\it Nonlinear Fokker-Planck Equations. Fundamentals and Applications}, Springer, Berlin. Heidelberg. New York, 2005.\vspace*{-3mm}


\bibitem{9a} Frank, T.D., Daffertshofer, A.,  $H$-theorem for nonlinear Fokker-Planck equations related to generalized thermostatics,{\it Physica A. Statistical Mechanics and its Applications}, 295(2001), 455-474.\vspace*{-3mm}
	

\bibitem{11b} Jordan, R., Kinderlehrer, D., Otto, F., The variational formulation of the Fokker-Planck equation, {\it SIAM J. Math. Anal.},  29 (1998), 1-17.\vspace*{-3mm}

\bibitem{18prim} Manita, O.A., Romanov, M.S., Shaposhnikov, S.V., On uniqueness of solutions to nonlinear Fokker-Planck-Kolmogorov equations, {\it Nonlin. Anal.}, 128 (2015), 199-226.\vspace*{-3mm}


\bibitem{18second} Manita, O.A., Shaposhnikov, S.V., Nonlinear parabolic  equations for measures,  {\it St. Petersburg Math. J.}, 25 (1) (2014), 43-62.\vspace*{-3mm}


\bibitem{16a} Markowich, P.A., Villani, C., On the trend to equilibrium for the Fokker-Planck equations: an interplay  between physics and functional analysis, {\it Mathematics Contemporary}, 2000.\vspace*{-3mm}

\bibitem{25a} Otto, F., Villani, C., Generalization of an inequality by Talagrand and links with the logarithmic Sobolev inequality, {\it J.~Funct. Anal.}, 173 (2000), 361-400.\vspace*{-3mm}

\bibitem{9} Pazy, A., The Lyapunov method for semigroups of nonlinear contractions in Banach spaces, {\it Journal d'Analyse Math.}, 40 (1981), 239-262.\vspace*{-3mm}

\bibitem{10} Schw\"ammle, V., Nobre, F.D., Curado, E.M.F., Consequences of the $H$-theorem from nonlinear Fokker-Planck equations, {\it Phys. Rev.}, E 76 (2007), 041123.\vspace*{-3mm}

\bibitem{26a} Temam, R., {\it Infinite Dimensional Dynamical System in Mechanics and Physics,} Springer-Verlag, New York. Berlin. Heidelberg. London. Paris. Tokyo, 1988.\vspace*{-3mm}

\bibitem{14prim} Trevisan, D., Well-posedness of multidimensional diffusion processes with weakly differentiable coefficients, {\it Electron. J.~of Probab.}, 21 (22) (2016), 1-41.\vspace*{-3mm}

\bibitem{22} Wang, F.-Y., {\it Functional inequalities, Markov semigroups and spectral theory},  Science Press, 2005, xx+379 pp. ISBN:7-03-014415-5.

\end{thebibliography}
\end{document}